\documentclass[ejs]{imsart}

\RequirePackage[OT1]{fontenc}
\RequirePackage{amsthm,amsmath}
\RequirePackage[numbers]{natbib}
\RequirePackage[colorlinks,citecolor=blue,urlcolor=blue]{hyperref}

\usepackage{graphicx}

\usepackage{amsfonts}
\usepackage{algorithm}
\usepackage{algorithmic}

\usepackage{latexsym}

\newtheorem{theorem}{Theorem}[section]

\newtheorem{lemma}{Lemma}[section]
\newtheorem{corollary}{Corollary}[section]
\newtheorem{remark}{Remark}[section]

\newtheorem{definition}{Definition}[section]

\newcommand{\R}{\mathbb{R}}

\newcommand{\PP} {{  \rm I\hskip-0.22em P}}

\newcommand{\EE} {{\rm I\hskip-0.48em E}}

\usepackage{natbib}

\pubyear{2014}
\volume{0}
\issue{0}
\firstpage{1}
\lastpage{8}

\startlocaldefs
\numberwithin{equation}{section}
\theoremstyle{plain}

\endlocaldefs

\begin{document}

\begin{frontmatter}
\title{On higher order isotropy conditions and lower bounds for sparse quadratic forms\thanksref{T1}}
\runtitle{Isotropy and quadratic forms}
\thankstext{T1}{The authors gratefully acknowledge financial support of the Swiss
National Science Foundation, grant nr.\ 20PA20E-134495. We are moreover very grateful to
Guillaume Lecu\'e for his helpful  comments.}

\begin{aug}
\author{\fnms{Sara} \snm{van de Geer}\ead[label=e1]{geer@stat.math.ethz.ch}}
\and
\author{\fnms{Alan} \snm{Muro}\ead[label=e2]{muro@stat.math.ethz.ch}}

\address{Seminar for Statistics\\
ETH Z\"urich\\
R\"amistrasse 101\\
8092 Z\"urich, Switzerland\\
\printead{e1,e2}}

%\author{\fnms{Third} \snm{Author}
%\ead[label=e3]{third@somewhere.com}
%\ead[label=u1,url]{www.foo.com}}
%
%\address{Address of the Third author\\
%usually few lines long\\
%usually few lines long\\
%\printead{e3}\\
%\printead{u1}}

%\thankstext{t1}{Some comment}
%\thankstext{t2}{First supporter of the project}
%\thankstext{t3}{Second supporter of the project}
\runauthor{S. van de Geer et al.}

\affiliation{Seminar for Statistics, ETH Z\"urich}

\end{aug}

\begin{abstract}
This study aims at contributing to lower bounds for 
empirical compatibility constants or empirical restricted
eigenvalues. This is of importance in compressed sensing and
theory for $\ell_1$-regularized estimators. Let $X$ be an $n \times p$ data matrix
with rows being independent copies of a $p$-dimensional random variable.
Let $\hat \Sigma := X^T X / n$ be the inner product matrix.
We show that the quadratic forms $u^T \hat \Sigma u$ are lower bounded by a value converging to one,
uniformly  over the set
of vectors $u$ with $u^T \Sigma_0 u $ equal to one and $\ell_1$-norm at most $M$.
Here $\Sigma_0 := \EE \hat \Sigma$ is the theoretical inner product matrix which we assume to exist. 
The constant $M$ is required to be of small order $\sqrt {n / \log p }$.
We assume moreover $m$-th order isotropy for some $m >2$ and sub-exponential
tails or moments up to order $\log p$  for the entries in $X$. 
As a consequence we obtain convergence of the empirical
compatibility constant to its theoretical counterpart, and
similarly for the empirical restricted eigenvalue. If the data matrix $X$
is first normalized so that its columns all have equal length we obtain lower bounds assuming only isotropy and
no further moment conditions on its entries.
The isotropy condition is shown to hold for
certain martingale situations. 
\end{abstract}

\begin{keyword}[class=MSC]
\kwd{62J07}
%\kwd{60K35}
%\kwd[; secondary ]{60K35}
\end{keyword}

\begin{keyword}
\kwd{compatibility}
\kwd{isotropy}
\kwd{quadratic forms}
\kwd{restricted eigenvalue}
\kwd{sparsity}
\kwd{transfer principle}
\end{keyword}
%\tableofcontents
\end{frontmatter}

\thispagestyle{plain}

\section{Introduction}\label{introduction.section}
Let $X$ be an $n \times p$ data matrix with rows being i.i.d.\ copies of a random
vector $X_0^T \in \R^p$.
We consider the empirical inner product matrix $\hat \Sigma= X^T X / n $.
For a vector $u \in \R^p$, let $\| u \|_q $ be its
$\ell_q$-norm ($1 \le q \le \infty $). We examine sparse quadratic forms
$u^T \hat \Sigma u$ where $u$ is sparse in the sense that $\| u \|_1 \le M$
for some constant $M \ge 1$.
We will provide lower bounds for 
$\min \{ u^T \hat \Sigma u : \ u^T \Sigma_0 u  = 1 , \ \| u \|_1 \le M \}  $
with $\Sigma_0 := \EE \hat \Sigma$  being the theoretical inner product matrix which we assume
to exist.
The constant $M$ will be required to be of small order
$\sqrt { n/ \log p }$. 
%In an asymptotic formulation we show in Theorem \ref{two-sided.theorem}
%that
%$$\hat \psi^2 (M) = \psi_0^2 (M) + o_{\PP} (1)  , $$
%where $\psi_0^2 (M)$ is the theoretical $\ell_1$-restricted eigenvalue
%$$\psi_0^2 (M) := \min \{ u^T \Sigma_0 u : \ \| u \|_2 = 1 , \ \| u \|_1 \le M \}  .$$

A motivation to study lower bounds for quadratic forms comes from theory for $\ell_1$-penalized
estimation methods. Here the so-called (empirical) compatibility constant plays an important role.
It is defined as follows.
For $S$ being a subset of $\{ 1 , \ldots  , p \}$, write $u_{j,S} = u_j {\rm l} \{ j \in S \}$
($j=1 , \ldots , p $) and $u_{-S} := u - u_S$.
The compatibility constant (\cite{vandeG07}) is
$$\hat \phi^2 (L,S):= \min \{ |S|  u^T \hat \Sigma u  :\  \| u_S \|_1 = 1 , \ \| u_{-S} \|_1 \le L \} . $$

The condition $\hat \phi^2 (L,S) >0$ for suitable values of $L$ and $S$ allows one to  establish oracle
inequalities for the Lasso. 
Indeed, let $u^0$ be the sparse vector we want to recover and let $S:= \{ j:\ u_j^0 \not= 0 \} $ be its active set.
Let $\xi \in \R^n$ be a ``noise" vector. Consider the Lasso 
$$\hat u := \arg \min_{u \in \R^p} \biggl \{ \| \xi + X u^0  - X u \|_{2,n}^2 + 2 \lambda \| u \|_1\biggr \} $$
where $\lambda>0$ is a tuning parameter and where we
use the notation $\| v \|_{2,n}^2 := v^T v / n $, $v \in \R^n$.
 For $\lambda > \lambda_0 =: \| \xi^T X \|_{\infty}/n $ it holds that
\begin{equation}\label{oracle.equation}
 \| X (\hat u -  u^0 ) \|_{2,n}^2 \le (\lambda + \lambda_0)^2  |S| / \hat \phi^2 (L,S) 
 \end{equation}
where $L:= (\lambda + \lambda_0)/(\lambda- \lambda_0 )$. We refer to \cite{BvdG2011} and the references therein. 
In the literature result (\ref{oracle.equation}) is considered to be an ``oracle inequality" if - in a suitable  asymptotic formulation - the constant $L$ remains bounded  (i.e. $\lambda$ is of the same order as $\lambda_0$) and 
$\hat \phi^2 (L,S)$ stays away from zero. In the present paper, this case
may serve as benchmark case. We give non-asymptotic results and
some asymptotic consequences showing that under certain conditions $\hat \phi^2 (L,S)$ indeed
stays away from zero.

Closely related is the so-called null space property\footnote{We thank Emmanuel Cand\'es for pointing this out.} (see e.g.\ \cite{gribonval2007highly}) used in exact recovery.
One says that $X$ has the null space property relative to $S$ if
for all $u \in \R^p$ with $Xu =0$ it holds that $\| u_S \|_1 < \| u_{-S} \|_1 $. 
The null space property is the same as  the
condition  $\hat \phi(1,S)>0$ and implies in the noiseless case
exact recovery of a sparse signal $u^0$ with active set $S$ using basis pursuit
(\cite{chen1998atomic}): 
$$\arg \min \{ \| u \|_1 : \ Xu = X u^0  \} = u^0 . $$
The compatibility constant  is also a close relative of the 
(empirical) restricted eigenvalue defined in \cite{bickel2009sal} as
$$\hat \kappa^2 (L, S) :=  \min \{   u^T \hat \Sigma u   :\ \| u \|_2 =1 , \ \| u_{-S} \|_1 \le L \| u_S \|_1  \} .$$
It is easy to see that $\hat \kappa^2  (L, S) \le \hat \phi^2 (L,S)$.
An example where $\hat \kappa^2  (L, S) $ will be much smaller than $ \hat \phi^2 (L,S)$ is
given in \cite{vdG:2009}. In that example $ \hat \kappa^2 (1,S)$ is about $1/|S|$ whereas
$\hat \phi^2 (1,S)$ is about $1/2$. One sees that for large values of $|S|$ the difference is substantial.

In many cases (e.g.\ when applying the Lasso) the data are first normalized: for $\hat \sigma_j^2 := \hat \Sigma_{j,j}$ one
replaces the $j$-th column $X_j$ of $X$ by $\tilde X_j/ \hat \sigma_j$, $j=1 , \ldots , p $. 
Therefore we study in Section \ref{normalized.section} the compatibility constant for normalized design
$$ \tilde \phi^2 (L,S) := \min \{ |S| u^T \hat R u : \ \| u_S \|_1 =1 , \ \| u_{-S} \|_1 \le L \} $$
and restricted eigenvalue for normalized design
$$\tilde \kappa^2 (L,S) := \min \{   u^T \hat R u   :\ \| u \|_2 =1 , \ \| u_{-S} \|_1 \le L \| u_S \|_1  \} .$$
where
$\hat R := \tilde X^T \tilde X / n $
with $\tilde X := ( \tilde X_1 , \ldots , \tilde X_p)$ being the normalized Gram matrix.

\subsection{Organization of the paper}\label{organisation.section}
After some notations and definitions in the next section, we present
in Section \ref{sparse.section} a bound for sparse quadratic forms.
The lower bounds for the empirical compatibility constant and
empirical restricted eigenvalue  follow from this. The upper bounds
depend on  fourth moments.
We will show that $\hat \phi (L,S)$ converges to
its theoretical counterpart, and similarly for $\hat \kappa (L,S)$ (see Theorem
\ref{two-sided.theorem}). For this we need $(L+1) \sqrt s$ to be of small
order $\sqrt {n / \log p}$ (for the lower bound). 
This is detailed in Section \ref{l1eigenvalue.section}. 
In Section \ref{normalized.section} we consider the transfer principle from \cite{oliveira2013} which
allows one to show that for the case where the data are normalized
very weak moment conditions suffice. Section \ref{related-work.section} is devoted to a discussion
with related work.  There we summarize the comparison of results in Table \ref{tab5}. 
In Section \ref{almost-bounded.section}
we make a brief comparison of the results when we drop
the isotropy assumption. 
We show convergence
of $| u^T \hat \Sigma u - u^T \Sigma_0 u |$ 
uniformly over $\| u \|_1 \le M $
assuming sub-exponential entries in $X_0$.
In Section \ref{martingale.section} we examine the higher order isotropy condition.
Finally, Section \ref{proofs.section} contains the proofs.

\section{Notation and definitions}\label{notation.section}

We let $\Sigma_0 := \EE X_0 X_0^T = \EE \hat \Sigma$ be the theoretical inner product matrix.
Its smallest eigenvalue is denoted by $\psi_0^2$. We do not assume $\psi_0 >0$.
For $m \ge 1$, and $Z$ a real-valued random variable,  we introduce the notation
$$ \| Z \|_{m}^m :=  \EE |Z|^m . $$
Thus $ u^T \Sigma_0 u = \| \bigl < X_0 u \bigr > \|_2^2 $ where $\bigl < X_0 u \bigr > $ is the inner
product $X_0^T u$, $u \in \R^p$. 

Let $X_{ i, \cdot}^T $ be the $i$-th row of $X$ ($i=1 , \ldots , n $). 
We write for a function $f : \ \R^p \rightarrow \R$,
$$\| f \|_{2,n}^2 := \| f (X) \|_{2,n}^2:=  {1 \over n}\sum_{i=1}^n f^2 (X_{i, {\cdot} })  $$
and
$$ \| f \|_2^2 := \| f (X_0) \|_{2}^2 := \EE  \| f (X) \|_{2,n}^2   $$
so that $\| X u \|_{2,n}^2 = u^T \hat \Sigma u $.
%We furthermore let for $k \ge 1$ and $\Psi_k (z) := \exp [ |z |^k ] -1$, $z \in \R$,
%the Orlicz-norm  $\| Z \|_{\Psi_k}$ be 
%$$\| Z \|_{\Psi_k} := \inf \{ c>0 :\  \EE \Psi_k ( | Z| / c ) \le 1 \} . $$

\begin{definition} \label{Bernstein.definition}
We say that a random variable $Z$ is Bernstein with constants $\sigma$ and $K$
if for all $k \in \{ 2,3, \ldots \}$
$$\EE | Z|^k \le { k! \over 2 } K^{k-2} \sigma^2 . $$
\end{definition} 

\begin{definition}\label{sub-Gaussian.definition}
We say that a random variable $Z \in \R$ is sub-Gaussian with constant $C$ if for
all $\lambda  >0 $
$$\EE \exp [ \lambda |Z| ] \le  2 \exp[ \lambda^2 C^2  / 2 ] . $$
\end{definition}
Let us denote for $k=1,2$ the Orlicz norm by
$$\| Z \|_{\Psi_k} := \inf\{c>0 :\  \EE \exp[|Z/c|^k] -1 \le 1 \} . $$
Then being a Bernstein random variable is equivalent to having finite $\| \cdot \|_{\Psi_1}$-norm (i.e., being
sub-exponential) and sub-Gaussianity is equivalent to a finite
$\| \cdot \|_{\Psi_2}$-norm. We have chosen for the Definitions \ref{Bernstein.definition}
and \ref{sub-Gaussian.definition} in order to have simple explicit dependence on the constants later on.

Note that if a random variable is sub-Gaussian with constant $C$ it is also
Bernstein with constants $\sigma = 2$ and $K=\sqrt 2C$.
Moreover, a Bernstein random variable $Z$ with constants $\sigma$ and $K$ always 
has $\sigma \le 3 K$ and so $\| Z \|_m \le m K $ for all $m \in \{ 3,4, \ldots \}  $. 

We use the definition of \cite{mendelson2008uniform} or \cite{loh2012} of a sub-Gaussian vector (in a slightly
alternative formulation). 
\begin{definition}\label{sub-Gaussianvector.definition}
A random vector $X_0 \in \R^p$ is sub-Gaussian with
constant $C$ if for all $u \in \R^p$ with $\| \bigl < X_0 u \bigr > \|_2= 1$ the random variable $\bigl < X_0 u \bigr >$ is sub-Gaussian with
constant $C$.
\end{definition}

The main concept we will use in this paper is weak isotropy for which we now present the definition.

\begin{definition}\label{weakisotropic.definition}
 Let $m \ge 2$. The  random vector $X_0 \in \R^p$ is weakly $m$-th order isotropic 
with constant $C_m$ if for all $u \in \R^p$ with $ \| \bigl < X_0 u \bigr > \|_2 =1$ it holds that
$$P ( | \bigl < X_0 u \bigr > | > t )  \le (C_m/ t)^{m}   \ \forall \ t >0. $$
\end{definition}

A Gaussian vector is sub-Gaussian with constant $1$ and is strongly
$m$-th order  isotropic (defined in Definition \ref{strongisotropic.definition}) with constant $\sqrt 2 \Gamma^{1 \over m} ((1+m)/2) \pi^{-{1 \over m} } $,
$ m \ge 2$. 

Definition \ref{weakisotropic.definition} (and \ref{sub-Gaussianvector.definition})
are invariant under rotations: 
if $\psi_0>0$ one may without loss of generality assume $\Sigma_0 =I$ here.
We however explicitly do not assume $\Sigma_0 =I$ because conditions on the $\ell_1$-norm
are not invariant under rotation. In contrast to the 
literature where the ``isotropic" case is sometimes defined as the case $\Sigma_0=I$
our definition of isotropy is rather to be understood as uniformity in all one-dimensional directions
(very much like isotropy of functions in Besov spaces).

\section{Lower bounds for sparse quadratic forms under higher order isotropy}\label{sparse.section}

The first result of Theorem \ref{easy.theorem} below is as in \cite{srivastava2013covariance} and is
given for completeness. It is only of interest when $p $ is smaller than $n$.
  The result is improved in  
 \cite{mendelson2013}. We refer to Section \ref{related-work.section} for a discussion.
 The second result of Theorem \ref{easy.theorem} extends the situation to the case where 
$p $ can be larger than $n$ but $\ell_1$-restrictions are invoked.
%The proof follows to some extend \cite{mendelson2013}. 
Here we need bounds on Rademacher averages.
A Rademacher sequence  is a sequence of 
independent random variables 
$\epsilon_1 , \ldots , \epsilon_n$ where each $\epsilon_i$ takes values $\pm 1$ with probability $1/2$.
We assume that $\epsilon:= (\epsilon_1 , \ldots , \epsilon_n )^T$ is
independent of $X$. Consider the Rademacher averages
$$W^T :=(W_1 , \ldots , W_p)^T:= \epsilon^T X / n . $$
and let $ \| W \|_{\infty}: = \max_{1 \le j \le p } | W_j |$.
We will need bounds for  $ \EE \| W \|_{\infty} $.
If the entries in $X_0$ are Bernstein with constants $\sigma_X$ and $K_X$, then
applying Lemma 14.12 in \cite{BvdG2011} gives
 \begin{equation}\label{expectationbound.equation}
 \EE \| W \|_{\infty}   \le 
 \sigma_X\sqrt {2\log (2p) \over n} + K_X { \log (2p) \over n }.
 \end{equation}
 and it is this bound that is invoked in the second result of Theorem \ref{easy.theorem}.
Further bounds for $\EE \| W \|_{\infty}$ are discussed in
Subsection \ref{EW.section}.

\begin{theorem} \label{easy.theorem}
Suppose that for some $m > 2$ the random vector $X_0$ is weakly $m$-th order isotropic with constant $C_m$
and define
\begin{equation}\label{Dm.equation}
D_m:= [ 2 C_m]^{m \over m-1} {(m-1)/( m-2)}  . 
\end{equation}
Then for all $t>0$ with probability at least $1- \exp[-t]$
\begin{equation}\label{first.equation}
    \inf_{\| \bigl < X_0 u \bigr > \|_2 =1  } \| X u \|_{2,n}^2 -1 \ge- \biggl [ 
     D_m
    \biggl ( 16 \sqrt {p \over n}     + \sqrt { 2 t \over n } \biggr )^{m-2 \over m-1}  
    + { 8 D_m^2 \over 3} \biggl ( {t \over n}  \biggr )^{m-2 \over m-1} \biggr ] .
    \end{equation}
    If in addition the entries in $X_0$ are Bernstein with constants $\sigma_X$ and $K_X$, then for all $t >0$ with probability at least $1- \exp [-t] $
\begin{equation}\label{second.equation}
    \inf_{\| \bigl < X_0 u \bigr > \|_2 =1 , \ \| u \|_1 \le M } \| X u \|_{2,n}^2 -1 \ge -\Delta_n^{\rm L} (M,t)
   \end{equation}
     where 
     \begin{equation} \label{Deltan.equation}
     \Delta_n^{\rm L}( M , t):=
     D_m
    \biggl ( 16 M \delta_n    + \sqrt { 2 t \over n } \biggr )^{m-2 \over m-1}  
    + { 8 D_m^2 \over 3} \biggl ( {t \over n}  \biggr )^{m-2 \over m-1} 
    \end{equation}
    with $\delta_n :=  \sigma_X\sqrt {2\log (2p) \over n} + K_X { \log (2p) \over n } $. 
    
            \end{theorem}
  
  \noindent          
 {\bf Asymptotics} In an asymptotic formulation suppose that $C_m$, $K_X$ and $\sigma_X$ remain fixed and
 that $(\log p)/n = o(1)$. Then the second result (\ref{second.equation}) of Theorem \ref{easy.theorem}
 says that for $M = o( \sqrt {n / \log p } )$ one has
$$  \inf_{\| \bigl < X_0 u \bigr > \|_2 =1 , \ \| u \|_1 \le M } \| X u \|_{2,n}^2 \ge 1 - o_{\PP} (1) . $$

 \begin{remark} The constant $t$ in the formulation of Theorem \ref{easy.theorem} allows one
 to choose the confidence level of the result. If $t$ is large (for example $p$ large and 
$ t= \log p $) the bounds will
 be true with large probability. Of course for very large $t$ the bounds will become void. 
 
 \end{remark}

% \begin{remark} \label{transformation.remark}
% In fact, Theorem \ref{easy.theorem} can be formulated in a more general fashion.
% The conditions on the entries of $X_0$ may be replaced by the same conditions on some
% non-degenerate linear transformation of $X_0$
% 
% \end{remark}

 \begin{remark} \label{constants.remark} We have not attempted to obtain small constants in the bound
 of Theorem \ref{easy.theorem}. In fact, the last ``smaller order" term in the expression 
 (\ref{Deltan.equation}) for $\Delta_n^{\rm L} (M,t)$ can be refined but this will make the expressions
 more involved. 
\end {remark}

% \begin{corollary} \label{p<n.corollary} As a special case one may consider the one
% where $p$ is at most $n$ 
% Then one finds a result as in \cite{srivastava2013covariance} (see also \cite{mendelson2013}
% for this result and for sharper bounds). Namely, 
% suppose that, for some $m > 2$, $X_0$ is weakly $m$-th order isotropic with constant $C_m$.
%Apply Theorem \ref{easy.theorem} with
%$D_m$ defined in (\ref{Dm.equation}) to get
%that
%for all $t >0$ with probability at least $1- \exp [-t] $
%$$
%    \inf_{\| \bigl < X_0 u \bigr > \|_2 =1 } \| X u \|_{2,n}^2   \ge 1 - \Delta_n^{\rm L} ( t) , $$
%    where 
%        $$\Delta_n^{\rm L} (t) := 
%    D_m
%    \biggl ( {16 }   \sqrt { p  \over n } + \sqrt { 2 t \over n } \biggr )^{m-2 \over m-1}  
%     +
%    { 8 D_m^{2 } \over 3} \biggl ( {t \over n}  \biggr )^{m-2 \over m-1} . $$
%   
% \end{corollary}
  \begin{remark} The technique to prove Theorem \ref{easy.theorem} does not rely on the fact that we
 consider squared functions $| (Xu)_i |^2$, $i=1 , \ldots , n$. For example, one may use it for bounding
 $$ \inf_{ \| \bigl < X_0 u \bigr > \|_2 =1 , \| u \|_1 \le M } \| X u \|_{q,n}^q  $$
 where for $q \ge 1$
  $$ \| X u \|_{q,n}^q := {1 \over n} \sum_{i=1}^n | ( Xu )_i |^q .  $$
  Then one could e.g.\ use weak isotropy conditions of order $m >q$.
  However a motivation for
  having such results is perhaps lacking.
 \end{remark}
 
 Theorem \ref{easy.theorem} is based on a truncation argument. For the case of a sub-Gaussian
 vector $X_0$ the truncation level can be taken rather small leading to an improved bound.
 We present this for completeness in the next lemma.
 
  \begin{lemma} \label{sub-Gaussianeasy.lemma}
  If the random vector $X_0$ is sub-Gaussian with constant $C$ we find
       that for all $t>0$ with probability at least $1- \exp[-t]$
  $$
  \inf_{\| \bigl < X_0 u \bigr > \|_2 =1 , \ \| u \|_1 \le M } \| X u \|_{2,n}^2 -1 \ge 
  $$ $$ - \sqrt 2 C b (1+ 2 \sqrt {2 \log(C/b) }) - 16 C^2 t \cdot { \log (C/b) \over 3n} $$
  where
  $$b:= 16 \min \biggl \{ M \delta_n^{\prime}  , \sqrt { p \over n } \biggr \} + \sqrt { 2 t \over n } . $$
 with $\delta_n^{\prime} :=  C \sqrt {2\log (2p) \over n} $.
 \end{lemma}

 \subsection{Bounds for $ \EE \| W \|_{\infty}$} \label{EW.section}
% If the entries in $X_0$ are Bernstein with constants $\sigma_X$ and $K_X$, then
%applying Lemma 14.12 in \cite{BvdG2011} gives
% \begin{equation}\label{expectationbound.equation}
% \EE \| W \|_{\infty}   \le 
% \sigma_X\sqrt {2\log (2p) \over n} + K_X { \log (2p) \over n }
% \end{equation}
%If these entries are sub-Gaussian with constant $C$ we find
%$$\EE \| W \|_{\infty}  \le C \sqrt {2\log (2p) \over n} . $$
% 
%
Inequality (\ref{expectationbound.equation}) presents a bound for $\EE \| W \|_{\infty}$
assuming Bernstein conditions. This bound is then invoked in Theorem \ref{easy.theorem}.
One may derive alternative bounds for $\EE \| W \|_{\infty} $ and adjust the definion of
$\delta_n$ in Theorem \ref{easy.theorem} accordingly. 
For example one may impose existence of $k$-th moments of
the entries of $X_0$ where $k$ is of order $\log p$. The paper
 \cite{lecue2014} presents
 refined results which we cite in the next lemma.
 
 \begin{lemma} \label{Lecue.lemma}
 Let $Z_1 , \ldots , Z_n$ be i.i.d. copies of a mean-zero random variable $Z \in \R$
 and $\bar Z := \sum_{i=1}^n Z_i / n $.
 Suppose that for some constants $\kappa_1$ and $\alpha \ge 1/2$ one has 
 $$ \| Z \|_k \le \kappa_1 k^{\alpha} ,\  2 \le k \le k_0 . $$
 Then for $n \ge k_0^{\max \{ 2 \alpha -1 , 1\} }$ and for all $k \le k_0$
$$ \| \bar Z \|_k \le c_0 \exp[2 \alpha -1 ]  \kappa_1 \sqrt {k/n}  $$
where $c_0$ is a universal constant.
 
 \end{lemma}
 
 \begin{corollary}\label{Lecue.corollary} Suppose that for some constants $\kappa_1$, $\eta\ge 2/ \log p $ and $\alpha \ge 1/2$ one has
 \begin{equation}\label{moments.equation}
\max_{1 \le j \le p } \| X_{0,j} \|_k \le \kappa_1  k^{\alpha} , \ 2 \le k \le k_0 := \eta \log p . 
\end{equation}
 Then for $n \ge k_0^{\max \{ 2 \alpha -1 , 1\} }$ 
 $$\max_{1 \le j \le p } \| W_j  \|_{k_0} \le c_0 \exp[2 \alpha -1 ]  \kappa_1 \sqrt {k_0/n} , $$
 where $c_0$ is a universal constant. But then
 $$ \EE \| W \|_{\infty} \le c_1 \sqrt { \log p / n }  $$
 where $c_1=
 c_0  \kappa_1 \sqrt {\eta} \exp [  2 \alpha -1 +1/ \eta ] $.

 \end{corollary}

 \section{Convergence of the compatibility constant and restricted eigenvalue}\label{l1eigenvalue.section}
 
% We begin with a remark which may be useful for the case where the
% data are normalized to have length one.
% 
%   \begin{remark}
%Let
% $$\hat \rho^2 (M) := \inf_{\| u  \|_2 =1 , \ \| u \|_1 \le M } u^T \hat R u $$
% where $\hat R:= \tilde X^T \tilde X  / n $ is the inner product matrix of the
% $\tilde X_{i,j} := X_{i,j} / \hat \sigma_j $ with $\hat \sigma_j^2 :=\hat \Sigma_{j,j}= \| X_j \|_{2,n}^2 $. 
% Assume for simplicity that  $\Sigma_{0, j,j} =1$ for all $j$.
% Then one can use Theorem \ref{easy.theorem} to find under its assumptions
% and with the notation used there that for all $t>0$ with probability
% at least $1- \exp[-t]$
% $$ \inf_{j} \hat \sigma_j^2 \ge 1- \Delta_n^{\rm L} (1, t)$$
%      Thus, with probability at least $1- \exp [-t]$ we have
%     $$\hat \rho^2 (M) \ge (1- \Delta_n^{\rm L} (1, t)) \hat \psi^2 (M) . $$
%     \end{remark}

 An ``almost isometric" (in a terminology from \cite{mendelson2013}) lower bound for the empirical 
 compatibility constant and empirical restricted eigenvalue follows easily from Theorem
 \ref{easy.theorem} as is shown in the next theorem.
 
 Recall that $S \subset \{ 1 , \ldots , p \}$ is an arbitrary subset. Let for $s:= |S|$
 $$ \phi_0^2 (L,S) := \min \{s  \| \bigl < X_0 u \bigr > \|_2^2:\ \| u_S \|_1 = 1 , \ \| u_{-S} \|_1 \le L  \}. $$
 be the theoretical compatibility constant and
 $$ \kappa_0^2 (L,S) := \min \{ \| \bigl < X_0 u \bigr > \|_2^2 : \ \| u_S \|_2 =1 , \ \| u_{-S} \|_1 \le L \| u_S \|_1 \} . $$
 be the theoretical restricted eigenvalue. 
 
 \begin{theorem} \label{psi.lemma}
  Under the conditions of Theorem \ref{easy.theorem} and using its notation we find that for all $t >0$,
 with probability at least $1- \exp[-t]$
 $$ 
   { \hat \phi^2(L,S) \over   \phi_0^2 (L,S) }-1  \ge-  \Delta_n^{\rm L} ( (L+1) \sqrt {s} / \phi_0 (L,S) , t)    $$
    as well as 
    $$ { \hat \kappa^2 (L,S)  \over \kappa_0^2 (L,S) } -1 \ge - \Delta_n^{\rm L} ( (L+1) \sqrt {s} / \kappa_0 (L,S),t )   . $$
      \end{theorem}

 Note that Theorem \ref{psi.lemma} does not depend on the smallest eigenvalue
 $\psi_0$ of $\Sigma_0$ nor on its maximal eigenvalue. If
 $\psi_0>0$ one may however want to insert the bounds $\phi_0 (L,S) \ge \kappa_0 (L,S) \ge \psi_0$.
 We refer to the ``Asymptotics" paragraph at the end of this section for a further discussion.

 %We denote the largest eigenvalue of $\Sigma_0$ by $\psi_{\rm max}^2$.
 % \begin{theorem} \label{isometric.theorem}
% Under the conditions of Theorem \ref{easy.theorem} and using its notation we find that for all $t >0$,
% with probability at least $1- \exp[-t]$
%\begin{equation} \label{isometricmain.equation}
%    \hat \psi^2 (M) \ge  \psi_0^2 (M)   -   \psi_{\rm max}  \Delta_n^{\rm L} (M/ \psi_0 ,t).
%    \end{equation}
%  \end{theorem}
% 
%  
%  The argument of Theorem \ref{isometric.theorem} cannot be used to show that
%  a lower bound for the compatibility constant  $\hat \phi^2(L,S)$ (restricted eigenvalue $\hat \kappa^2 (L,S)$) 
%is asymptotically close to its  theoretical counterpart. This is due to the fact that
%  the restriction $\| u_{S^c } \|_1 \le L \| u_S \|_1 $ does not lead to any
%  good upper bound for $\| u_{S^c } \|_2 $. (We used in Lemma \ref{compatibility.lemma} that
% $$\{ \| u_{S^c } \|_1 \le L  \| u_S \|_1 =\sqrt s   \} \subset 
%\{ \| u_{S^c } \|_1 \le L  \| u_S \|_1= \sqrt s , \ \| u_S \|_2 \ge 1   \}  $$ $$ \subset \{ \| u \|_1 \le (L+1) \sqrt {s} , \| u \|_2 \ge 1   \} . $$ The infimum over the latter set is attained in for values $\| u \|_2=1$.) 
%In other words, the price for having a good upper bound
%for $\| u_{S^c} \|_2$ is the gap between $\hat \phi(L,S)$ (and $\hat \kappa (L,S)$)
%and $\hat \psi^2 (M)$,

  The next 
 issue is whether 
$\hat \phi(L,S)$ actually converges to $\phi_0(L,S)$ and $\hat \kappa (L,S)$
to $\kappa_0(L,S)$. This part follows easily from the lower bounds of Theorem \ref{psi.lemma} and convergence
of 
 $\| X u \|_{2,n}^2 - \| \bigl < X_0 u \bigr > \|_2^2$ for fixed values of $u$, for which in turn
  we e.g.\ would like to have fourth moments.  If $m > 4$, this $4$-th order
 moment condition follows from $m$-th order weak isotropy. If however $X_0$
 is only $m$-th order weakly isotropic for $m\le 4$ we need some other
 means to check $4$-th moments. The next lemma can be invoked.
   
 \begin{lemma}\label{fourthmoments.lemma}   Suppose that  the entries in $X_{0}$ are Bernstein with constants $\sigma_X$ and $K_X \ge \sigma_X$ and
 that for some constant $c_0\ge 1$ and for $c_1= 2(1+c_0)(K_X + \sigma_X)$ we have
 $$c_1 M  \log (2p) \le p^{c_0/2}   . $$
 Then for all $u$ with $\| \bigl < X_0 u \bigr > \|_2 \le 1$ and $\| u \|_1 \le M$ we have
 $$ \| \bigl < X_0 u \bigr > \|_4^2 \le \sqrt 2 c_1 M \log (2p)    . $$
 \end{lemma}

Combining Theorem \ref{psi.lemma} with Lemma \ref{fourthmoments.lemma}  gives
the upper and lower bounds shown in the next theorem.
 
 \begin{theorem} \label{two-sided.theorem} Suppose that
$X_0$ is weakly 
 $m$-th
 order isotropic with constant $C_m$ and that the entries in $X_0$ are Bernstein
 with constants $\sigma_X$ and $K_X > \sigma_X$.
 For the case
 $m\le 4$ we assume in addition that 
 for some constant $c_0\ge 1$ and for $c_1:= 2(1+c_0) (K_X + \sigma_X)$
 \begin{equation}\label{c1.equation}
  c_1 (L+1) \sqrt s   \log (2p) \le p^{c_0/2}    .
  \end{equation}
Define $D_m$ as
 in (\ref{Dm.equation}) and $\Delta_n^{\rm L}(M,t) $ as in (\ref{Deltan.equation}).
   For all $t >0$, with probability at least $1- \exp [-t] - 1/t$
$$ - \Delta_n^{\rm L} ( (L +1) \sqrt s  / \phi_0 (L,S) , t)    \le { \hat \phi^2 (L,S) \over
\phi_0^2 (L,S) }  -  1  \le 
  \Delta_n^{\rm U} ((L +1) \sqrt s ,t) $$
and
$$ - \Delta_n^{\rm L} ( (L+1) \sqrt s  / \kappa_0 (L,S) , t)    \le { \hat \kappa^2 (L,S) \over  \kappa_0^2 (L,S)} - 1  \le 
  \Delta_n^{\rm U} ((L +1)  \sqrt s ,t) $$
where
$$\Delta_n^{\rm U} ((L+1) \sqrt s ,t) = \begin{cases}   c_1 (L+1) \sqrt s  \log (2p)\sqrt {2t \over n } 
 &  m \le 4 \cr
\min\biggl  \{  c_1 (L+1) \sqrt s   \log (2p)    , C_m^2 \sqrt { {m \over 2(m-4)} }   \biggr  \}  \sqrt {2t \over n } & m > 4 \cr 
 \end{cases}.
  $$
 
   \end{theorem}

 {\bf Asymptotics} In an asymptotic formulation we assume that the constants $1/ \phi_0(L,S)$, $C_m$, $\sigma_X$ and
 $K_X$ remain bounded. Then it follows from Theorem \ref{two-sided.theorem} that
 under its conditions, 
 as long as $(L+1) \sqrt s = o(  \sqrt {n/ \log p  }  ) $ 
  $$\hat \phi^2 (L, S) \ge  \phi_0^2 (L,S)- o_{\PP }(1) . $$
If $m \le 4$ , $c_0$ is fixed and
$(L+1) \sqrt s= o(  \sqrt n / \log p )$ we find
$$\hat \phi^2 (L,S) \le \phi_0^2 (L,S) + o_{\PP }(1) . $$

Similar results hold for the restricted eigenvalue.
(Note that for $c_0$ fixed  and $p > n$ condition (\ref{c1.equation}) follows from the already
imposed condition
$(L+1) \sqrt s = o(  \sqrt {n/ \log p  }  ) $.)
Thus, in the upper bound an additional $\sqrt {\log p }$ appears in the requirement on $M$. 
This term can be omitted 
if $m >4$ or if we assume the entries in $X_0$ are sub-Gaussian instead of Bernstein. 

%  Under  the sparsity assumption $\| u^0 \|_1 \le M $ where $u^0 := \arg \min \{ \| \bigl < X_0 u \bigr > \|_2: \ \|u \|_2=1  \} $
% one of course has $\psi_0^2 (M) = \psi_0^2$, i.e. then we can conclude that
% $\hat \psi^2 (M) $converges to $\psi_0^2$. 
% 
%

 \section{Bounds for the  compatibility constant and restricted eigenvalue using the transfer principle}
\label{normalized.section}

In this section, we assume for simplicity that $\Sigma_0$ has ones on the diagonal.
We let $\hat \sigma_j^2 := \hat \Sigma_{j,j}= \| X_j \|_{2,n}^2$, $j=1, \ldots , p$ where
$X_j$ denotes the $j$-th column of $X$. 

\subsection{The transfer principle}
The transfer principle given in the next theorem is from \cite{oliveira2013}. As shown in the latter paper
it can be used to move from the case $p\le n$ to $p > n$ assuming $\ell_1$-conditions.
We will apply this technique here as well, for non-normalized design in Theorem
\ref{transfer1.theorem} and for normalized design in Theorem \ref{normalized.theorem}.
The results are compared with \cite{oliveira2013} in Section \ref{related-work.section}.

\begin{theorem} \label{transfer.theorem}
Let $A$ be a symmetric $p \times p$ matrix with $A_{j,j} \ge 0 $ for all $j \in \{1 , \ldots , p \}$. 
Let $d \in \{ 2 , \ldots , p \}$  and suppose that for all
$J \subset \{ 1 , \ldots , p \}$ with cardinality $|J| = d $ and all $u \in \R^p$ one has
$$u_J^T A u_J \ge 0 . $$ Then for all $u \in \R^p$
$$ u^T A u \ge - \max_{j} A_{j,j}  \| u \|_1^2 / (d-1) . $$
\end{theorem}

We will invoke the transfer principle via the following corollary (as well as directly in the proof
of Theorem \ref{normalized.theorem}).
The corollary  is as in \cite{oliveira2013} and we state it here in our notation for ease of reference.

\begin{corollary} \label{transfer.corollary}
Let $M^2 \in \{ 2 , \ldots , p \} $ and $0 \le \Delta < 1 $.
Consider the events
$${\cal A}:= \biggl \{ {  \| X u_J \|_{2,n}^2   \over \| \bigl < X_0 u \bigr >_J \|_2^2 }  \ge 1 - \Delta   ,\
\forall \ u \in \R^p, \  \forall \ J \subset \{ 1 , \ldots , p \} \ {\it with} \  \ | J | \le M^2  \biggr \}  $$
and, for some $\epsilon >0$, the event
$$ {\cal B}:= \biggl \{ \max_{1 \le j \le p} \hat \sigma_j^2  \le 1+ \epsilon  \biggr  \} . $$
Then on ${\cal A}\cap {\cal B}$ for all $\| u \|_1 \le M \| \bigl < X_0 u \bigr > \|_2 $
$$ \| X u \|_{2,n}^2 \ge (1- 3 \Delta - 2 \epsilon ) \| \bigl < X_0 u \bigr > \|_2^2  . $$
 \end{corollary}
 
 To put this corollary to work we insert 
the first result of Theorem \ref{easy.theorem}. 
 
  \begin{theorem}\label{transfer1.theorem}
 Suppose that for some $m > 2$ the random vector $X_0$ is weakly $m$-th order isotropic with constant $C_m$.
Define $D_m$ as in Theorem \ref{easy.theorem}.
Let for $M^2 \in \{2, , \ldots , p \}$
$$
 \bar \Delta_n^{\rm L}( M, t):=
     D_m
    \biggl (  16M  { (1+ \sqrt {2 \log p}) \over \sqrt n }+ \sqrt { 2 t\over n } \biggr )^{m-2 \over m-1}  
    + { 8 D_m^2 \over 3} \biggl ( {t + M^2 \log p  \over n}  \biggr )^{m-2 \over m-1} $$
    and let, for some $\epsilon >0$, ${\cal B}$ be the event
    $${\cal B} := \{ \max_{1 \le j \le p } \hat \sigma_j^2  \le1+ \epsilon \} . $$
    
Then with probability at least $1-  \exp [-t] - \PP ({\cal B})$ uniformly in $\| u \|_1 \le M \| \bigl < X_0 u \bigr > \|_2 $
$$  \| X u \|_{2,n}^2 \ge \biggl (1- 3 \bar \Delta_n^{\rm L} (M, t) - 2\epsilon \biggr  )  \| \bigl < X_0 u \bigr >\|_2^2 .$$
                 \end{theorem} 
The above theorem invokes  Theorem \ref{easy.theorem}  for handling
 the event ${\cal A}$. One may also use            
 the results in \cite{oliveira2013} for the case of
 $m$-th order strong isotropy (defined in  Definition \ref{strongisotropic.definition}) {with $m \ge 4$ or those which can be deduced from  \cite{mendelson2013} for
 the case $m$-th order weak isotropy with $m >2$ (the latter paper does not explicitly
 treat an event of the form ${\cal A}$). 
 For the case $ m < 4$ for instance the
arguments in \cite{mendelson2013}
would allow  to replace $\bar \Delta_n (M,t)$ in Theorem \ref{transfer1.theorem}
 (which is of order $\bigl [ M\sqrt {\log p /n}\bigr ]^{m-2 \over m-1}$ by
 a term of order $\bigl [  M  \sqrt {\log p / n } \log ( 1/ (\sqrt {\log p / n } M))  \bigr ]^{2(m-2) / m}$.

Clearly, one can again apply the results to the compatibility constant and restricted eigenvalue
as in Theorem \ref{psi.lemma}. This gives
the following corollary.

\begin{corollary}Assume the conditions of Theorem \ref{transfer1.theorem} and
let $(L +1)^2 s \in \{ 2 , \ldots , p \}$. Then with probability at least $1- \exp[-t] - \PP ({\cal B})$
$$ { \hat \phi^2 (L,S) \over \phi_0^2 (L,S)} \ge 1- 3 \bar \Delta_n^{\rm L} ((L+1) \sqrt s/ \phi_0 (L,S) , t) -2 \epsilon $$
as well as
$$ { \hat \kappa^2 (L,S) \over \kappa_0^2 (L,S)} \ge 1- 3 \bar \Delta_n^{\rm L} ((L+1) \sqrt s/\kappa_0 (L,S)  , t) -2 \epsilon .$$
\end{corollary} 

\subsection{The behaviour of $\max_j \hat \sigma_j^2 $}\label{sigma.section}

Recall that in Theorem \ref{easy.theorem} the lower bound for
$\inf \{ \|  X u \|_{2,n}: \ \| \bigl < X_0 u \bigr > \|_2 =1 , \ \| u \|_1 \le M \} $ depends on
the bound 
$\delta_n $ for $\EE \| W \|_{\infty} $. Bounding $\EE \| W \|_{\infty}$ leads to moment conditions
on the entries in $X_0$. 
The transfer principle now
leads to requiring a bound for $\max_j \hat \sigma_j^2 $ where
$\hat \sigma_j^2 = \sum_{i=1}^n X_{i,j}^2/n $. The latter is clearly a more
difficult task than the former.  In the non-normalized case this appears 
 to be the price to pay for
application of the elegant transfer principle.

We first assume sub-Gaussian tail behaviour in Lemma \ref{sigmaGauss.lemma} and
then moments up to order $\log p$ in Lemma \ref{sigmamoments.lemma}. 
 
 \begin{lemma} \label{sigmaGauss.lemma}Suppose that the entries of $X_0$ are sub-Gaussian with constant $C$.
 Then for all $t >0$
 $$ \PP \biggl ( \max_{1 \le j \le p} | \hat \sigma_j^2 - 1 |/ (8C^2)  \ge
 \sqrt {2 \log (2p) \over  n }  + \sqrt {2t \over n}+ { \log (2p) \over n }     + {t\over n } 
\biggr ) \le \exp [-t ] . $$
 
 \end{lemma}

 \begin{lemma}\label{sigmamoments.lemma} Suppose the conditions of Corollary \ref{Lecue.corollary}
with constants $\kappa_1\ge 1 $, $\eta  \ge 2/ \log p $ and $\alpha \ge 1/4$:
 $$\max_{1 \le j \le p } \| X_{0,j} \|_k \le \kappa_1  k^{\alpha} , \ 2 \le k \le k_0 := \eta \log p . $$
 Then for $n \ge (k_0/2)^{\max \{ 4 \alpha -1 , 1\} }$ and all $t >0$ with probability at least $1-1/t$
 $$\max_{1 \le j \le p }\hat \sigma_j^2 \le 1 +  c_1 t^{2 / k_0} \sqrt { \log p  / n }  $$
 where 
 $ c_1:= c_0 \exp[ 4 \alpha -1+ 2/\eta ] 
 \kappa_1^2 2^{2 \alpha+1}\sqrt {\eta/2}$
 with $c_0$ a universal constant.  

\end{lemma}

  For example, when $\alpha=1$, the direct approach of
  Theorem \ref{easy.theorem} requires $n \ge \log p $ (see Corollary \ref{Lecue.corollary})
  whereas the transfer
  principle leads to requiring $n \ge \log^3§  p$.

\subsection{Normalized design}

Define  $\tilde X_j := X_j / \hat \sigma_j , \ j=1 , \ldots , p$,
$\tilde X := ( \tilde X_1 , \ldots , \tilde X_p)$  and
$$ \hat R := \tilde X^T \tilde X / n .   $$
Thus $ \hat R= \hat D^{-1/2} \hat \Sigma \hat D^{-1/2}$ 
where $\hat D= {\rm diag} ( \hat \sigma_1^2 , \ldots , \hat \sigma_p^2 ) $. 
Define for $S \subset \{ 1 , \ldots , p \}$ a set with cardinality $s := |S|$ the 
(empirical) compatibility constant for normalized design
$$ \tilde \phi^2 (L,S) := \min \{ s u^T \hat R u :\ \| u_S \|_1 =1 , \ \| u_{-S} \|_1 \le L \} . $$
Similarly, the (empirical) restricted eigenvalue for normalized design is
$$\tilde \kappa^2 (L,S) := \min \{ u^T \hat R u :\ \| u_S \|_2=1 , \ \| u_{-S} \|_1 \le L \| u_S \|_1 \} . $$

In  \cite{BvdG2011} the (theoretical) adaptive restricted eigenvalue is defined as
$$\kappa_*^2 (L,S) := \min \{ \| \bigl < X_0 u \bigr > \|_2^2: \ \| u_S \|_2 =1 , \ \| u_{-S}\|_1 \le L \sqrt s \} .  $$
Clearly $\kappa_*^2 ( L, S) \le \kappa_0^2(L,S)$.
We prove in Theorem \ref{normalized.theorem} that the empirical compatibility constant $\tilde \phi^2(L,S)$ can be
bounded from below by the theoretical adaptive restricted eigenvalue. The theorem establishes that
 compatibility needs no further moment conditions on the entries
in $X_0$.
If we do assume such  moment conditions on the entries $X_{0,j}$ with $j \in S$, the results
can be extended to restricted eigenvalues, as shown in \cite{oliveira2013} for the case of $4$-th order
strong isotropy (defined in Definition \ref{strongisotropic.definition}), and as shown in the next theorem.

\begin{theorem}\label{normalized.theorem}
Suppose that for some $m > 2$ the random vector $X_0$ is weakly $m$-th order isotropic with constant $C_m$.
Define $D_m$ as in Theorem \ref{easy.theorem}, and let $\bar \Delta_n^{\rm L} ( M, t)$ be defined
as in Theorem \ref{transfer1.theorem}. 
Let, for some $\epsilon >0$ sufficiently small, ${\cal B}_S$ and ${\cal C}_S $ be the events
\begin{equation}\label{CS.equation}
{\cal B}_S := \biggl \{ \max_{j \in S} \hat \sigma_j^2 \le 1+\epsilon \biggr \} , \ 
{\cal C}_S:= \biggl \{ \sum_{j \in S}  \hat \sigma_j^2 /s \le 1+\epsilon\biggr \} . 
\end{equation}
Let $0<\Delta <1$ be arbitrary and define $M^2(\Delta)$ as the largest value
of $M^2 \in \{ 2 , \ldots , p\}$ such that $\bar \Delta_n^{\rm L} ( M(\Delta)  , t) \le \Delta$
(assuming such a value exists).  Define $L(\Delta, \epsilon) := L \sqrt {1+ \epsilon} / (1- \sqrt {\Delta} )$.
Then with probability at least $1- \exp[-t] - \PP ({\cal C}_S ) $ we have
$$ \tilde \phi^2 (L,S) \ge { \kappa_*^2 (  L(\Delta, \epsilon) , S)  \over  (1+ \epsilon) } -
{ (L+1)^2 s \over M^2 (\Delta)-1 } . $$
Moreover, with probability at least $1- \exp[-t] - \PP({\cal B}_S)$
 $$ \tilde \phi^2(L,S) \ge { \phi_0^2 (  L(\Delta, \epsilon) , S)  \over  (1+ \epsilon) } -
{ (L+1)^2 s \over M^2 (\Delta)-1 } $$
and
$$ \tilde \kappa^2 (L,S) \ge { \kappa_0^2 (  L(\Delta, \epsilon) , S)  \over  (1+ \epsilon) } -
{ (L+1)^2 s \over M^2(\Delta)-1 } . $$
%Here, $M_2 (\Delta)$ is defined as $M_1 (\Delta)$ but with $\phi_0 (L,S)$ replaced
%by $\kappa_0 (L,S)$.
%

\end{theorem}

{\bf Asymptotics} The above theorem shows that when $1/ \kappa_* (L,S) = {\mathcal O}(1)$,
$C_m = {\mathcal O}(1)$ and 
$(L+1) \sqrt s = o ( \sqrt {n/\log p })$
then also $1/ \tilde \phi (L,S) = {\mathcal O}_{\PP} (1)$
and in fact $\liminf_{n \rightarrow \infty} \tilde \phi(L,S) / \kappa_* (L,S) \ge 1$
since $\kappa_* (L, S)$ is continuous in $L$
and $\sum_{j \in S} \hat \sigma_j^2 / s = 1+ o_{\PP} (n^{-{m-2 \over 2}} )$.

\section{Related work} \label{related-work.section}
Before discussing related work we present the definitions of the concepts used.
Recall that in this paper we require weak isotropy (see Definition \ref{weakisotropic.definition}).

\begin{definition}\label{strongisotropic.definition}
 Let $m \ge 2$. The  random vector $X_0 \in \R^p$ is strongly $m$-th order isotropic 
with constant $C_m$ if for all $u \in \R^p$ with $ \| \bigl < X_0 u \bigr > \|_2 =1$ it holds that
$$\| \bigl < X_0 u \bigr > \|_m \le C_m . $$
\end{definition}

\begin{definition}\label{L1L2.condition}
The  random vector $X_0 \in \R^p$ satisfies the $L_1$-$L_2$ property
with constant $C$ if for all $u \in \R^p$ with $\| \bigl < X_0 u \bigr > \|_2=1$ it holds that
$$\| \bigl < X_0 u \bigr > \|_1 \ge 1/C . $$
\end{definition}

\begin{definition}\label{smallball.definition}
The  random vector $X_0 \in \R^p$  satisfies the small ball property with constants
$C_1>0$ and $C_2>0$ if for all
$u \in \R^p$ with $\| \bigl < X_0 u \bigr > \|_2=1$ it holds that
$$ P\bigl  ( | \bigl < X_0 u \bigr > | \ge 1/ C_1 \bigr) \ge 1/ C_2  . $$
\end{definition}

It can be shown that for appropriate constants one has (for $m > 2$)
$${\mbox {strong $m$-th order isotropy}} \Rightarrow {\mbox {weak $m$-th order isotropy} }  $$ $$
\Rightarrow  {\mbox {$L_1$-$L_2$ property}} \Rightarrow {\mbox{small ball property} }.$$
E.g. for the last implication see \cite{mendelson2013}.

\subsection{Relation of this work with \cite{mendelson2013} and  \cite{oliveira2013}}
The paper \cite{mendelson2013} obtains lower bounds for the smallest eigenvalue of
$\hat \Sigma$ for the case $p \le n$. Their approach allows one to
show that for $p \ll n$ it holds that $u^T \hat \Sigma u \ge (1- \Delta) u^T \Sigma_0 u$
uniformly in $u \in \R^p$
with large probability for some small $\Delta$. Such a result is not stated explicitly but
it is easy to infer.
The bounds in 
 \cite{mendelson2013} are better than the first result (\ref{first.equation}) of Theorem \ref{easy.theorem}.
 The paper employs a type of ``peeling device" and the fact that
 for all $0<a<b$
 $$\{ \{ x:\ | x u | > K \} :\ \| \bigl < X_0 u \bigr > \|_2 =1 , \ K  \in (a, b]  \}  $$
 is a VC-class with dimension at most $p$. If we have
``good" bounds for the entropy for $\| \cdot \|_{2,n}$ of the classes
 $$\{ \{ x:\ | x u | > K \} :\ \| \bigl < X_0 u \bigr > \|_2 =1 , \ \| u \|_1 \le M , \ K  \in (a, b]  \}  $$
their argument can be extended to the case $p>n$
 with $\ell_1$-restrictions.
However,  how to derive ``good" entropy bounds for such classes is as yet not clear to us.
 
 Both papers  \cite{mendelson2013} and \cite{oliveira2013}
assume $m$-th order isotropy (defined here in Definitions \ref{weakisotropic.definition} and 
\ref{strongisotropic.definition}).  The paper \cite{mendelson2013} has results with weak isotropy for any
$m>2$, whereas \cite{oliveira2013} assumes strong isotropy with $m=4$. 
The paper \cite{oliveira2013} shows that by a
transfer principle (described here in Theorem \ref{transfer.theorem}) a result for $p\le n$ can be
invoked to derive that also for the case $p \gg n$ one has
$u^T \hat \Sigma u \ge (1- \Delta^{\prime} ) u^T \Sigma_0 u $ uniformly
in $\|u \|_1 \le M u^T \Sigma_0 u$ with large probability for some small $\Delta^{\prime}$
and not too large $M$ (generally of small order $\sqrt {n / \log p }$). 
%\begin{remark}
%One can partly improve the bound (\ref{second.equation})
% of Theorem \ref{easy.theorem} by applying
% the results of \cite{mendelson2013} for the case $p \le n$ and then apply the transfer principle
% of \cite{oliveira2013}. However, we then find
% that we need a stronger moment conditions on the entries of $X_0$, see Section \ref{normalized.section}. 
%  \end{remark}
 In the present paper we consider weak isotropy with $m>2$ as in \cite{mendelson2013} and we 
show by a direct method that $u^T \hat \Sigma u \ge (1- \Delta) u^T \Sigma_0 u $ 
uniformly in $\| u \|_1 \le M u^T \Sigma_0 u $ with large probability for some
small $\Delta$. Here, we assume sub-exponential tails
for the entries in $X_0$, or, inserting results from \cite{lecue2014}, existence of moments
up to order $\log p$ for these entries. We compared the result with the
one using the transfer principle of \cite{oliveira2013}. Our finding is that the
transfer principle needs slightly stronger moment conditions. In fact, our direct approach requires
a bound for the maximum of the $p$ Rademacher averages of the columns of $X$,
whereas the approach using the transfer principle makes it necessary
to have a bound for the maximal length of the $p$ columns of $X$. Both can be dealt with
by assuming higher order moments, but clearly the Rademacher averages need
less moments than the lengths. 

The paper \cite{oliveira2013}  shows that when
the columns of $X$ are normalized to have all equal length, then
the transfer principle leads to lower bounds for the (empirical) compatibility
constant and (empirical)
restricted eigenvalues assuming only $4$-th order strong isotropy and moments of order bigger than 4
for the entries in $X_0$. We presented this result in Section \ref{normalized.section} relaxing
$4$-th order strong isotropy to $m$-th order weak isotropy with $m>2$. Moreover,
we derive that the compatibility constant $\hat \phi^2 (L,S)$ is positive with
large probability assuming only isotropy but no additional moment assumptions on the entries
in $X_0$. Thus, using normalized design one obtains exact recovery under
isotropy only.

\subsection{Further related work}
In \cite{rudelson2011reconstruction}
a result of \cite{raskutti2010restricted} concerning a lower bound for restricted eigenvalues is extended from the Gaussian case
to the sub-Gaussian case. The paper \cite{adamczak2010quantitative} considers the case
of log-concave distributions, which is related to sub-exponentiality of the vector $X_0$
(the sub-exponential variant of Definition \ref{sub-Gaussianvector.definition}). 
The papers \cite{srivastava2013covariance}
and \cite{mendelson2013} provide lower bounds
for the empirical smallest eigenvalue $\hat \psi^2  := \min \{ u^T \hat \Sigma u : \ \| u \|_2 =1 \} $ for the case where $p $ is at most
$n$. The paper \cite{srivastava2013covariance} uses higher order isotropy conditions (defined in Definitions \ref{strongisotropic.definition} and \ref{weakisotropic.definition}) and the paper \cite{mendelson2013} uses these too, but they in addition explore
small ball properties (defined in Definition \ref{smallball.definition}).
 The paper \cite{lecue2014} considers the null space property and
 restricted eigenvalues $\hat \kappa^2 (L,S)$ invoking small ball properties.
 Indeed, they show that small ball properties are very natural requirements when one aims
 at lower bounds.
 With the small ball property one obtains  an ``isomorphic" bound (we call this a result of type II in Table \ref{tab5}), that is,
 in a standard asymptotic framework the lower bound remains
 strictly smaller than the theoretical counterpart. 
 Apart from the small ball property the paper  \cite{lecue2014} needs moment conditions. It requires
  the stronger (``sub-Gaussian") 
 conditions of Lemma \ref{sigmamoments.lemma} instead of the (``sub-exponential") condition
 (\ref{moments.equation}) of Corollary \ref{Lecue.corollary}. 
 The papers \cite{Lecuenecessary} and \cite{lecue2014} show that moment conditions are necessary
 for exact recovery. 
%The paper \cite{oliveira2013} includes results for compatibility constants and
%restricted eigenvalues when the data are first normalized. In that paper $4$-th order isotropy
%is assumed, but the moment conditions on the entries in $X_0$ are reduced to moments of
%order bigger than $4$ only. 

In Table \ref{tab5} we present a summary of the results in the cited papers in comparison with the present paper.
Of course it is not possible to make a simple comparison doing
all aspects of the cited papers justice. The summary should be seen as focussing on what are in our view
the relevant differences. 

\begin{table}[!htb]
% latex table generated in R 3.0.2 by xtable 1.7-1 package
% Tue Oct  1 13:12:03 2013
\begin{center} 
%Active set $S_0 = \{1,2,3\}$
\begin{tabular}{llllll}
  \hline
 & \cite{oliveira2013}&\cite{mendelson2013}&\cite{lecue2014}&  pp   \\ 
 \hline
 isotropy & $m \ge 4$ & $m >2$ & no & $m > 2$  \\ 
small ball & no & yes & yes & no  \\ 
  conditions on $p$ & no  & $p < n$ & no & no  \\ 
  sub-Gaussian & yes & yes & no & no  \\ 
   normalized & yes & no & no & yes  \\ 
   results &  $\kappa^2$ & $\psi^2$  & $\phi^2$ & $\kappa^2$ and $\phi^2$  \\
   type of result &  I & I \& II &   II & I  \\ 
   moment conditions: &  &  &  &   \\ 
   \ \ non-normalized  & - & -  & sub-Gaussian type& sub-exponential type  \\ 
   \ \ normalized  & weak for $\kappa^2$ & - & - & none for $\phi^2$  \\
  \hline
\end{tabular}
\end{center}
\caption{The entry ``pp" stands for the present paper. With ``isotropy" we mean
weak or strong isotropy. 
The ``sub-Gaussian" results  concern the lower tails for quadratic forms.
With ``conditions on $p$" we mean conditions stronger than the
asymptotic one $\log p / n \rightarrow 0$. 
The entries ``normalized" stand for normalized design and ``non-normalized" for
non-normalized design. 
The symbols $\kappa^2$, $\psi^2$ and $\phi^2$ are shorthand for restricted eigenvalue,
smallest eigenvalue and
compatibility constant respectively. With results of ``type I" we mean
results in terms of theoretical counterparts. Results of ``type II" are in terms of the
constants occurring e.g.\ in the small ball property.
The ``moment conditions" are apart from isotropy (or small ball properties) those on the entries of $X_0$. 
%The entry ``necessary" in the ``pp" column merely means that we point out in the present paper that we can use the same
%moment conditions as in \cite{lecue2014} which are proved to be necessary in \cite{Lecuenecessary.}
}\label{tab5}
\end{table}

 \section{The case of (almost) bounded random variables}\label{almost-bounded.section}

 The bounded case is considered \cite{rudelson2011reconstruction} and
 a reformulation is in 
 \cite{vandeGeer2014}.
 It is shown there  that when $\| X_0 \|_{\infty} \le K_X$ then for a universal
 constant $c_1$ and  for all $t>0$, with probability at least $1- \exp[-t] $
 $$ \sup_{\| \bigl < X_0 u \bigr > \|_2 \le 1 , \ \| u \|_1 \le M } 
 \biggl | \| X u \|_{2,n} - \| \bigl < X_0 u \bigr > \|_2^2 \biggr | / c_1  \le M K_X \sqrt { \log p \log^3 n + t \over n }
 $$ $$+ M^2 K_X^2 { \log p \log^3 n + t \over n } . $$
 Observe this inequality goes both ways, and it does not require
 higher order isotropy conditions. On the other hand, the bounds involve an additional
 $\log^3 n$-factor.  If we replace the assumption
 of bounded random variables by (say) a sub-Gaussian assumption
 but {\it do} assume strong (say) isotropy we can again use
 a truncation argument and
 obtain an inequality that goes both ways. Admittedly, the number of $\log p$- and $\log n$-terms
 increases. 
 
 We first present an auxiliary truncation lemma. 
\begin{lemma}\label{truncation.lemma} 
Suppose $X_0$ is strongly $m$-th order isotropic with constant $\tilde C_m$ and that its components
are sub-Gaussian with constant $C$. Let $t >0$ be arbitrary and
let 
$$A (t):= \{x \in \R^p: \  \max_{1 \le j \le p } | x_j | \le C (\sqrt {2t+2  \log (2 p)+ 2 m(\log n ) / (m-2)} ) \} . $$
Then for all $u$ with $\| \bigl < X_0 u \bigr > \|_2=1$
$$\| (\bigl < X_0 u \bigr >) {\rm l} \{ X_0 \notin A (t)  \}  \|_2^2 \le \tilde C_m^2  \exp[- t (m-2) / m ] /n .$$
\end{lemma}

\begin{theorem} \label{uniform.theorem}
Suppose $X_0$ is strongly $m$-th order isotropic with constant $\tilde C_m$ and that its components
are sub-Gaussian with constant $C$. Then for a universal constant $c_1$ and for 
all $t>0$ with probability at least $1- (1+ n^{-m/(m-2)} )  \exp[-t]  $
 $$ \sup_{\| \bigl < X_0 u \bigr > \|_2 \le 1 , \ \| u \|_1 \le M } 
 \biggl | \| X u   \|_{2,n}^2 - \| \bigl < X_0 u \bigr >  \|_2^2 \biggr | / c_1  $$ $$\le M C
  \sqrt { ( 2t +2 \log (2p) + 2m (\log n)/ (m-2) ) (\log p \log^3 n + t)  \over n }
 $$ $$+ M^2  C^2 { ( 2t +2 \log (2p) + 2m (\log n) / (m-2) )  ( \log p \log^3 n + t ) \over n } $$ 
 $$+ \tilde C_m^2  \exp[- t(m-2)/m] /n  . $$

\end{theorem}

\section{Higher order isotropy}\label{martingale.section}

If $X_0$ is (strongly or weakly) $m$-th order isotropic with constant $C_m$ and $A$ is a 
$q \times p$ matrix, then clearly
$A X_0$ is also (strongly or weakly) $m$-th order isotropic with constant $C_m$.
In other words, the property is invariant under linear transformations.
The same is true for sub-Gaussianity. 
In particular, we have invariance under any permutation of the $X_{0,j}$.

In the next subsection, we assume that the $\{ X_{0,j} \}$ 
form  a directed acyclic graph (possibly after some linear transformation)
where the noise terms are a martingale difference array with fixed sub-Gaussian
tail behaviour. Then we extend in Subsections \ref{sub-Gauss.section} and
\ref{Bernstein.section} the situation where the 
conditional tail behaviour is sub-Gaussian or Bernstein, with constants
depending on predictable random variables. 
We
consider there a filtration $\{{\cal F}_j \}_{j=1}^p $ and predictable random variables
$\{ V_j \}_{j=1}^p $ that satisfy for some constants $m > 2$ and $\mu_m$
$$ \max_{1 \le j \le p } \| V_j \|_m \le \mu_m  . $$
We investigate strong $m$-th order isotropy. In fact we give explicit
expressions for $\| \bigl < X_0 u \bigr > \|_m$ in terms of 
$\| u \|_2 $. This implies strong isotropy if we assume the smallest eigenvalue
$\psi_0^2$ of $\Sigma_0$ is positive.
Obviously this also implies  a bound for the largest eigenvalue $\psi_{\rm max}^2$ of
$\Sigma_0$:
$$\psi_{\rm max}^2 \le \max\{ \| \bigl < X_0 u \bigr > \|_m^2: \ \| u \|_2 =1 \} . $$

\subsection{Directed acyclic graphs}\label{DAG.section}

Let $X_0$ be a vector of random variables with mean zero and covariance matrix
$\Sigma_0 := \EE X_0 X_0^T $. We want to find conditions
such that for all $u \in \R^p$, with $\| \bigl < X_0 u \bigr >  \|_2 =1$ the random variable $\bigl < X_0 u \bigr >$ is sub-Gaussian with
constant $C$.
We will examine this here for the situation where the graph of $X_0$ has a directed acyclic graph (DAG) structure 
that is, satisfying (after an appropriate permutation of the indexes) the structural equations model
\begin{equation} \label{SEM.equation}
X_{0,1} = \epsilon_{0,1} , \ X_{0,j}= \sum_{k=1}^{j-1} X_{0,k}  \beta_{k,j}  + \epsilon_{0,j} , \ j=2 , \ldots , p 
\end{equation}
where $\{ \epsilon_{0,j} \}_{j=1}^p $ is a martingale difference array for the filtration $\{ {\cal F}_j \}_{j=0}^{p-1} $.
We assume $X_{0,j}$ is ${\cal F}_j$-measurable, $j=1 , \ldots , p$.
We moreover assume that $\omega_j^2 := {\rm var} ( \epsilon_{0,j} )=
\EE {\rm var} (X_j \vert {\cal F}_{j-1}) $ exists for all $j$.
Note that model (\ref{SEM.equation}) holds when $X_0$ is Gaussian for example.
More generally, 
the standard linear structural equations model is a special case.
The latter model assumes that for $j \ge 2$, the noise $\epsilon_{0,j}$ is independent of
$\{X_{0,k} \}_{k=1}^{j-1}$, and that $\epsilon_{0,1} , \ldots \epsilon_{0,p}$ are independent
mean-zero random variables.

\begin{lemma} \label{DAG.lemma}Assume the structural equations model (\ref{SEM.equation}). Assume in addition
that for some constant $C$ and for all $\lambda \in \R$
$$ \EE ( \exp[ \lambda \epsilon_{0,j}/ \omega_j ]  \vert {\cal F}_{j-1} ) \le \exp [ \lambda^2 C^2 / 2 ] , 
j=1 , \ldots  , p . $$
Then $X_0$ is sub-Gaussian
with constant $C   $.
\end{lemma}

The above lemma follows from the fact that its condition implies that the vector $\epsilon_0:= ( \epsilon_{0,1} , \ldots , \epsilon_{0,p} )^T$ is
sub-Gaussian with constant $C$. If $\epsilon_0$ is (strongly or weakly) $m$-th order
isotropic with constant $C_m$, then under the structural equations  model (\ref{SEM.equation})
the vector $X_0$ is also  (strongly or weakly) $m$-th order
isotropic with constant $C_m$. This follows from the fact that $X_0$ is a linear transformation of
$\epsilon_0$. One may use the results of the next two subsections to check isotropy
of $\epsilon_0$.

\subsection{The conditionally sub-Gaussian case}\label{sub-Gauss.section}

Let $\{ {\cal F}_j \}_{j=0}^p$ be a filtration and for $j=1 , \ldots , p$, let
$X_{0,j}$ be ${\cal F}_j$-measurable and $V_j$ be ${\cal F}_{j-1}$-measurable.
We assume that for some $m >2$,
$$ \max_{1 \le j \le p } \| V_j \|_m := \mu_m < \infty . $$

\begin{lemma}  \label{sub-Gauss.lemma} Suppose that for all $j$ 
$$\EE (X_{0,j} \vert {\cal F}_{j-1} ) =0 , \ \EE ( \exp [ \lambda X_{0,j} ] \vert {\cal F}_{j-1} ) \le \exp [ \lambda^2 V_j^2 / 2 ] \  \forall \ \lambda \in \R . $$
If $\{ V_j \}_{j=1}^p $ is ${\cal F}_0$-measurable then for all $\| u \|_2=1$
$$ \| \bigl < X_0 u \bigr > \|_m \le \sqrt {2m}  \mu_m . $$
For general predictable $\{ V_j \}$ we have for $2 < m_0 < m $ and all
$ \| u \|_2 =1$
$$ \| \bigl < X_0 u \bigr > \|_m \le \sqrt {2 m \over m- m_0}   \biggr ( { 3 m \Gamma ( m_0 /2 +1 ) \over
m- m_0 } \biggr )^{1 / m_0}   \mu_m . $$
\end{lemma}

\subsection{The conditionally Bernstein (or sub-exponential) case}\label{Bernstein.section}
Let as in the previous sub-section $\{ {\cal F}_j \}_{j=0}^p$ be a filtration and for $j=1 , \ldots , p$, let
$X_{0,j}$ be ${\cal F}_j$-measurable and $V_j$ be ${\cal F}_{j-1}$-measurable and
satisfying  for some $m >2$,
$$ \max_{1 \le j \le p } \| V_j \|_m := \mu_m < \infty . $$
As in the previous section, 
we prove strong isotropy but now under a different condition. 

\begin{lemma} \label{Bernstein.lemma} Suppose that for some constant $K$ and all $j$
$$ \EE (X_{0,j} \vert {\cal F}_{j-1} ) = 0 , \ \EE( | X_{0,j} |^k \vert {\cal F}_{j-1}) \le { k! \over 2} K^{k-2} V_j^2 , \ k=2,3, \ldots . $$
If the $ \{ V_j \}_{j=1}^p$ are non-random, then for all $\| u \|_2 =1$
$$\|  \bigl < X_0 u \bigr > \|_m  \le  \sqrt {2m} \mu_m   +   m K . $$
If $\{ V_j \}_{j=1}^p $ is ${\cal F}_0$-measurable we get for all $\| u \|_2 =1$
$$ \| \bigl < X_0 u \bigr > \|_m \le 2^{1-1/m_0} \sqrt {2m} \mu_m   +   2^{1-1/m_0} m  K. $$
For general predictable $\{ V_j \}_{j=1}^p$ we have for all $2 < m_0 < m$
and all $\|u \|_2=1$
$$  \|  \bigl < X_0 u \bigr >  \|_{m_0} \le 
  \sqrt { 2m \over m- m_0 }  \biggl ( { 3 m \Gamma (m_0 /2 + 1)  \over m-m_0}  \biggr )^{m_0 /2 +1}  \mu_m+ 
\biggr (3  \Gamma (m_0 + 1) \biggr )^{1/ m_0} K . $$

\end{lemma}

Note that the conditions of the above lemma imply that the entries in $X_0$ are Bernstein
with constants $ \mu_2 $ and $K$, where $\mu_2 := \max_{1 \le j \le p } \| V_j \|_2 \le \mu_m$.
In other words, the conditions of the lemma imply the bound of
Theorem \ref{easy.theorem} with $\delta_n = \mu_2 \sqrt {2 \log (2p)/n} + K \log (2p)/n $ and with
$m$ replaced by any $m_0<m$.

\section{Proofs} \label{proofs.section}

%\subsection{Proof of the inequalities in Section \ref{introduction.section}}
%
%{\bf Proof of Lemma \ref{compatibility.lemma}.}
%Write $s:= |S|$. We apply the inequalities $\| u_S \|_1 \le \sqrt s \| u_S \|_2 \le
%\sqrt s \| u \|_2 $. Then we find
%$$\hat \phi^2 (L,S) = \min \{ s u^T \hat \Sigma u / \| u_S \|_1^2 : \| u_{-S} \|_1 \le L \| u_S \|_1 \} $$
%$$\ge \min \{ u^T \hat \Sigma u / \| u_S \|_2^2 :  \| u_{-S} \|_1 \le L \| u_S \|_1 \} 
% = \hat \kappa^2 (L,S)  $$ 
% $$ = \min \{ u^T \hat \Sigma u : \ | u_{-S} \|_1 \le L \| u_S \|_1 , \ \| u_S \|_2 = 1 \} $$
% $$\ge \min \{ u^T \hat \Sigma u  :  \| u \|_1 \le (L+1) \sqrt {s} , \| u_S \|_2 =1   \} $$
%$$\ge \min \{ u^T \hat \Sigma u  :  \| u \|_1 \le( L+1) \sqrt {s} , \ \| u \|_2 \ge 1   \} $$
%$$ =  \min \{ u^T \hat \Sigma u  :  \| u \|_1 \le( L+1) \sqrt {s}, \  \| u \|_2 = 1   \} 
%= \hat \psi^2 ( \sqrt s (L+1)) .  $$
%\hfill $\sqcup \mkern -12mu \sqcap$

\subsection{Proofs for Section \ref{sparse.section}}

Recall that Theorem \ref{easy.theorem} presents  lower bounds for sparse quadratic forms.

       {\bf Proof of Theorem \ref{easy.theorem}.}  
For $Z \in \R$, and $K >0$, we introduce the truncated version
$$[Z]_K := \begin{cases} -K , & Z < -K  \cr \ Z, & |Z|\le K \cr +K , & Z > K \cr  \end{cases} . $$

We obviously have for any $K >0$ and $u \in \R^p$
\begin{equation} \label{truncation.equality}
\| X u \|_{2,n}^2 \ge \| [ Xu]_K \|_{2,n}^2 
\end{equation} 
where $[Xu ]_K $ is the vector $\{ [ (Xu)_i]_K: \ i=1 , \ldots , n \} $ with,
for $i\in \{ 1 , \ldots , n \}$, $(Xu)_i$ be the $i$-th component of the vector $Xu$. 
Moreover, whenever $\| \bigl < X_0 u \bigr > \|_2 = 1$ by the weak isotropy
$$1 - \| [\bigl < X_0u \bigr >  ]_K \|_2^2 \le   2 C_m^m K^{-(m-2)}  / (m-2) . $$
Here, we used the formula
$$1 - \| [\bigl < X_0u \bigr >  ]_K \|_2^2=
\int_0^{\infty} \PP ( | \bigl < X_0u \bigr >  | > \sqrt { K^2 + t } ) dt . $$
We note that 
$$ \EE \sup_{\| \bigl < X_0 u \bigr > \|_2 =1  , \ \| u \|_1 \le M } \biggl | \| [Xu]_K \|_{2, n}^2 -  \| [\bigl < X_0u \bigr > ]_K \|_2^2 \biggr | $$
$$ ={1 \over K^2} \EE \sup_{\| \bigl < X_0 u \bigr > \|_2 =1/K  , \ \| u \|_1 \le M/K } \biggl | \| [Xu]_1 \|_{2, n}^2 -  \| [\bigl < X_0 u \bigr >]_1 \|_2^2 \biggr | .$$
 
 Let
 $${\bf Z} :=  \sup_{\| \bigl < X_0 u \bigr > \|_2 =1/K  , \ \| u\|_1 \le M/K } \biggl | \| [Xu]_1 \|_{2,n}^2 -  \| [\bigl < X_0u
 \bigr > ]_1 \|_2^2 \biggr | . $$

 By symmetrization (see e.g.\ \cite{vanderVaart:96}, p.108) and contraction (\cite{Ledoux:91}, p.112), 
 $$ \EE{\bf Z} 
  \le 2 \EE \sup_{ \| \bigl < X_0 u \bigr > \|_2 =1/K , \ \| u \|_1 \le M/K } \biggl | {1 \over n} \sum_{i=1}^n \epsilon_i [(Xu)_i ]_1^2 \biggr | $$
  $$ \le 8  \EE \sup_{\| \bigl < X_0 u \bigr > \|_2 = 1/K , \  \| u \|_1 \le M/K } \biggl | {1 \over n} \sum_{i=1}^n \epsilon_i (Xu )_{i} \biggr | $$
  since the mapping $Z \mapsto [Z]_1^2 $ is $2$-Lipschitz.  
Continuing with the last bound, we will apply 
 $$  \EE \sup_{ \| \bigl < X_0 u \bigr > \|_2 = 1/K  } \biggl | {1 \over n} \sum_{i=1}^n \epsilon_i (Xu)_{i} \biggr |
  ={1 \over K} \sqrt {p \over n}   $$
 for deriving (\ref{first.equation}) 
and
$$  \EE \sup_{   \| u \|_1 \le M/K } \biggl | {1 \over n} \sum_{i=1}^n \epsilon_i (Xu)_{i} \biggr | =
 {M \over K} \EE \| W \|_{\infty} \le { M \over K} \delta_n 
 $$
for deriving (\ref{second.equation}). 
  In other words
 $$ \EE {\bf Z} \le  {8 \over K}  \min \biggl \{ M  \delta_n ,   \sqrt {p \over n } \biggr \}. $$
 
   Next we apply the concentration inequality of \cite{Bousquet:02} to ${\bf Z}$.
   We get for all $t >0$
   $$\PP \biggl (  {\bf Z} \ge \EE {\bf Z} + { 2 t \over 3n} + \sqrt {2 t/n } \sqrt {1/K^2+ 4 \EE {\bf Z} }  \biggr ) \le \exp [-t] $$
   where we used for $\| \bigl < X_0 u \bigr > \|_2 \le 1/K$ the bound
   $${\rm var} ( [ \bigl < X_0 u \bigr >]_1^2 ) \le \EE [\bigl < X_0 u \bigr > ]_1^2  \le  \| \bigl < X_0 u \bigr > \|_2^2 \le 1/K^2  . $$
   We invoke that
    $$ \sqrt {2 t/n } \sqrt {1/K^2+ 4  \EE {\bf Z} }\le \sqrt {2 t/n }(  1/K + 2  \sqrt{   \EE {\bf Z} })$$
    $$ \le {\sqrt {2 t/n } \over K } + {2 t \over n} + \EE {\bf Z} . $$
   This gives for all $t >0$
   $$\PP \biggl (  {\bf Z} \ge 2 \EE {\bf Z}  +{ \sqrt {2 t/n }  \over K} + { 8 t \over 3n}  \biggr ) \le \exp [-t] $$
   and hence
   $$\PP \biggl (  {\bf Z} \ge {16 \over K} \min \biggl \{ M \delta_n ,  \sqrt { p \over n } \biggr \} 
   +{ \sqrt {2 t/n }  \over K} + { 8 t \over 3n}   \biggr ) \le \exp [-t] .
   $$ 
   So with probability at least $1- \exp[-t]$
   $$\inf_{\| \bigl < X_0 u \bigr > \|_2 =1 , \ \| u \|_1 \le M } \| X u \|_{2,n}^2-1  \ge  - {2 C_m^m \over (m-2) K^{m-2} }  
   $$ $$- 16 K \min \biggl \{ M \delta_n , 
  \sqrt { p \over n }  \biggr \} -K { \sqrt {2 t \over n } } - { 8 K^2 t \over 3n}.  $$
   We now let
   $$K:= [2 C_m^m ] ^{1 \over m-1} b^{-{1 \over m-1} } , $$
   where
   $$b:= 16 \min \biggl \{ M \delta_n  , \sqrt { p \over n } \biggr \} + \sqrt { 2 t \over n } . $$
   Then
   $${2 C_m^m  \over (m-2 )K^{m-2} } + K b = D_m b^{m-2 \over m-1}, $$
   and
   $${ 8 K^2 t \over 3n } \le  { 8 D_m^{2 } \over 3} \biggl ( {t \over n}  \biggr )^{m-2 \over m-1} . $$
   It follows that with probability at least $1- \exp [-t] $
    $$\inf_{\| \bigl < X_0 u \bigr > \|_2 =1 , \ \| u \|_1 \le M } \| X u \|_{2,n}^2 -1 $$ $$\ge - D_m
    \biggl ( 16 \min \{    M \delta_n , \sqrt { p \over n }\biggr \} + \sqrt { 2 t \over n } \biggr )^{m-2 \over m-1}  -
     { 8 D_m^{2 } \over 3} \biggl ( {t \over n}  \biggr )^{m-2 \over m-1}  .$$
     
 \hfill $\sqcup \mkern -12mu \sqcap$
 
  \begin{remark}
  With assumptions weaker than the weak isotropy assumption used in the present paper, for
example with the $L_1$-$L_2$ property, one can prove
 lower bounds  along the same lines as for Theorem \ref{easy.theorem}.
One applies instead of the truncation
 inequality 
 (\ref{truncation.equality}) in the proof of Theorem \ref{easy.theorem} the inequality
 $$ \| X u \|_{2,n} \ge \| Xu \|_{1,n} , $$
 where 
 $$ \| X u \|_{1,n} = {1 \over n} \sum_{i=1}^n | ( Xu )_i |  $$
 One can then proceed using the arguments following (\ref{truncation.equality}) in the proof of Theorem \ref{easy.theorem}
 using the Lipschitz property of the absolute
 value function $Z \mapsto |Z| $. For results assuming only the small ball property, we refer
 to  \cite{lecue2014}.
  \end{remark}

 We now provide a proof for the sub-Gaussian case along the same lines as the proof of
 Theorem \ref{easy.theorem}.
  
 {\bf Proof of Lemma \ref{sub-Gaussianeasy.lemma}.} 
 We use the same notation as in the proof of Theorem \ref{easy.theorem} for truncation at a value
 $K$. 
    Whenever $\| \bigl < X_0 u \bigr > \|_2 = 1$,
  $$ 1- \| [\bigl < X_0 u \bigr > ]_K \|_2^2 \int_0^{\infty}=
  \PP ( |\bigl < X_0 u \bigr > | > \sqrt { K^2 + t } ) dt  $$
  $$ \le 2 \int_{0}^{\infty} \exp[ - (K^2 +t ) / (2 C^2 ) ] = 4 C^2 \exp [ -K^2 / (2 C^2 ) ] . $$
  We choose
  $$ K =  C \sqrt { 2 \log (C/b) } $$
  where
  $$b:= 16 \min \biggl \{ M \delta_n^{\prime}  , \sqrt { p \over n } \biggr \} + \sqrt { 2 t \over n } . $$
The result then follows by the same arguments as those used for Theorem \ref{easy.theorem} and
  inserting that in the sub-Gaussian case one has $ \EE \| W \|_{\infty} \le \delta_n^{\prime} $.
  \hfill $\sqcup \mkern -12mu \sqcap$

  \subsection{Proofs for Section \ref{l1eigenvalue.section}}
  
   We first proof the
  ``almost isometric" bound for the compatibility constant and restricted eigenvalue.
  
        {\bf Proof of Theorem \ref{psi.lemma}.}
  By Theorem \ref{easy.theorem} we know that uniformly in $u$ with $\| u \|_1 \le
  M \| \bigl < X_0 u \bigr > \|_2 $ with probability at least $1- \exp[-t]$
  $$ \| X u \|_{2,n}^2 \ge ( 1- \Delta_n^{\rm L} ( M  , t ) ) \| \bigl < X_0 u \bigr >\|_2^2 . $$
  If $\| u_S \|_1 =1 $ and $\| u_{-S}  \|_1 \le L$ we clearly have
  $$\| u \|_1 \le (L +1) = (L+1) \| u_S \|_1 \le  (L+1) \sqrt s   \| \bigl < X_0 u \bigr > \|_2/ \phi_0 (L,S).$$
  This implies the lower bound for the compatibility constant.
  If $\| u_S\|_2 =1$ and $\| u_{-S} \|_1 \le L \| u_S \|_1 $ we again have
  $\| u \|_1 \le (L+1) \sqrt {s} \| u \|_2  \le (L+1) \sqrt s  \| \bigl < X_0 u \bigr > \|_2/ \kappa_0 (L,S) $ which implies
  the result for the restricted eigenvalue. 
 
  \hfill $\sqcup \mkern -12mu \sqcap$

%  Next is the ``almost isometric" bound for the $\ell_1$-restricted eigenvalue.
%  
%   {\bf Proof of Theorem \ref{isometric.theorem}.}
% Using the same truncation notation as in the proof of Theorem
% \ref{easy.theorem}, we have for all $K >0$
% $$ \hat \psi^2 (M) \ge
% \inf_{\| u \|_2 = 1 , \ \| u \|_1 \le M } \| ( Xu )_K \|_{2,n}^2 $$
% $$ = \inf_{\| u \|_2 = 1 , \ \| u \|_1 \le M } \biggl \{ \| \bigl < X_0 u \bigr > \|_2^2 + \biggl ( \| (\bigl < X_0 u \bigr >_K \|_2^2 - \| \bigl < X_0 u \bigr > \|_2^2 \biggr ) $$ $$
%\ \ \ \ \ \ \ \ \ \ \ \ \ \ \ \ \ \ \ \ \ \ + \biggl ( \| ( Xu )_K \|_{2,n}^2 - \| (\bigl < X_0 u \bigr >)_K \|_2^2 \biggr ) \biggr\} $$
%$$ \ge \psi_0^2 (M) - \sup_{\| u \|_2 = 1 , \ \| u \|_1 \le M } \biggl | \| (\bigl < X_0 u \bigr >)_K \|_2^2 - \| \bigl < X_0 u \bigr > \|_2^2 \biggr |
%$$ $$ - \sup_{\| u \|_2 = 1 , \ \| u \|_1 \le M } \biggl | \| ( Xu )_K \|_{2,n}^2 - \| (\bigl < X_0 u \bigr >)_K \|_2^2 \biggr | . $$
%If $\| u \|_2 = 1$, we have $\| \bigl < X_0 u \bigr > \|_2 \le \psi_{\rm max} $. Hence then, for
%$\tilde u := u/ \| \bigl < X_0 u \bigr > \|_2 $
%$$ \| \bigl < X_0 u \bigr > \|_m = \| X_0 \tilde u \|_m \| \bigl < X_0 u \bigr > \|_2 \le \psi_{\rm max} C_m . $$
%The theorem therefore follows by the same arguments as those for
%Theorem \ref{easy.theorem}.
%
%\hfill $\sqcup \mkern -12mu \sqcap$

We now check the fourth moments, i.e. the second moments of quadratic forms. 

 {\bf Proof of Lemma \ref{fourthmoments.lemma}.}
 
% This follows from the fact that for real-valued random variables $Z_1 , \ldots , Z_p$,
% one has
% $$\| \max_{1 \le j \le p } Z_j \|_m \le c_m \| \max_{1 \le j \le p } Z_j  \|_{\Psi_1} \le
% c_m^{\prime} \log p \max_{1 \le j \le p } \| Z_j \|_{\psi_1} $$
% where $c_m$ and $c_m^{\prime}$ depend only on $m>0$.
% We present a full proof for completeness.
 
  One readily sees that each $X_{0,j}$ has Orlizc norm $\| \cdot \|_{\Psi_1} $ bounded
 by $K_X + \sigma_X$:
 $$ \EE \exp [ | X_{0,j} | / ( K_X + \sigma_X ) ] -1 \le 1 ,\ \forall \ j . $$
 Hence for all $t>0$ and all $j $
 $$ \PP ( | X_{0,j} | > t (K_X + \sigma_X)) \le 2 \exp[ - t ] .$$
 It follows that for all $t >0$
 $$ \PP \biggl ( \max_{j  } |X_{0,j} | > [ t + (1+c_0) \log (2p) ] (K_X + \sigma_X) \biggr ) \le \exp[-t] / (2p)^{c_0}. $$
  Clearly 
 $$ \| \bigl < X_0 u \bigr >\|_4^4 = \underbrace{\EE |\bigl < X_0 u \bigr > |^4 {\rm l} \{ \max_{j  } |X_{0,j} | \le  2(1+c_0) \log (2p)  (K_X + \sigma_X) \} }_{:= i } $$
 $$ + \underbrace{\EE |\bigl < X_0 u \bigr > |^4 {\rm l} \{ \max_{j   } |X_{0,j} | >  2(1+c_0) \log (2p)  (K_X + \sigma_X) \} }_{:= ii}.$$
 We have for $\| u \|_1 \le M $ and $\| \bigl < X_0 u \bigr > \|_2 \le 1 $
 $$ i \le (2(1+c_0)M)^2 (K_X + \sigma_X )^2 \log^2 (2p) \EE \| \bigl < X_0 u \bigr > \|_2^2 $$ $$\le (2(1+c_0) M)^2 (K_X + \sigma_X )^2  \log^2 (2p) . $$
 Now for a random variable $Z$ satisfying for all $t>0$ 
 $ \PP ( |Z| > bt + K/2) \le c \exp[-t]$  for certain constants $b$, $c$ and  $K$
 $$ \EE Z^4 {\rm l} \{ |Z |> K \} = \EE ( Z^4 - K^4 ) {\rm l}\{ |Z| > K \} +
 K^4 \PP ( |Z| > K ) $$
 $$ = \int_0^{\infty} \PP ( Z^4   >t  + K^4 ) dt + K^4 \PP ( |Z| > K ) $$
 $$ \le \int_0^{\infty} \PP ( Z  >(t/8)^{1/4}  + K/2 ) dt  + K^4 \PP ( |Z| > K ) .$$ 
 Here we used that
 $(t^{1/4} + K)^4 \le 8( t + K^4) $ and $8^{1/4} \le 2 $. 
 So we get
 $$ \EE Z^4 {\rm l} \{ |Z |> K \} \le 
 8 b^4 \int_0^{\infty} \PP ( Z  >bs  + K/2 ) d s^4  + K^4 \PP ( |Z| > K )$$
 $$ \le 8b^4 c \int_{0}^{\infty} \exp[-s] d s^4 + c K^4 = (8 \times 4!) b^4 c + cK^4\exp[-K/(2b)] 
$$ $$ \le (4b)^4 c + c K^4 . $$
 Apply this to $|Z|:= \max_{1 \le j \le p } |X_{0,j} | $. Then we can take $b= (K_X + \sigma_X)$, 
 $c= 1/(2p)^{c_0} $ and $K= 2(1+ c_0)\log (2p)  (K_X + \sigma_X) $.
 We find
 $$ \EE \max_{1 \le j \le p } | X_{0,j }|^4 {\rm l } \{ \max_{1 \le j \le p } |X_{0,j} | >  2(1+c_0) \log (2p)  (K_X + \sigma_X) \} $$
 $$
 \le [4(K_X + \sigma_X)]^4 /(2p)^{c_0} + \biggl (2(1+c_0) \log (2p)  (K_X + \sigma_X)\biggr )^4 /(2p)^{c_0} $$ $$\le
 \biggl (2(1+c_0) \log (2p)  (K_X + \sigma_X)\biggr )^4/p^{c_0} 
  $$
  since $\log (2p) \ge1 $ and hence $2(1+c_0)\log (2p) \ge 4 $.
 But then
 $$ ii \le M^4 \EE \biggl ( \max_{1 \le j \le p } | X_{0,j} |^4 \{ \max_{1 \le j \le p } |X_{0,j} | >  2(1+c_0) \log (2p)  (K_X + \sigma_X) \}   \biggr )  $$
 $$ \le [2(1+c_0)M]^4  (K_X + \sigma_X)^4 \log^4 (2p)/p^{c_0}  $$ $$\le
 [2(1+c_0) M]^2 (K_X + \sigma_X)^2 \log^2(2p)   $$
 where in the last step we invoked the assumption of the lemma.
 We conclude
 $$ \| \bigl < X_0 u \bigr > \|_4^4 \le i + ii \le 2 [2(1+c_0) M]^2 (K_X + \sigma_X)^2 \log^2(2p) . $$
 \hfill $\sqcup \mkern -12mu \sqcap$

  As a result, we can now obtain lower and upper bounds for the compatibility constant and 
restricted eigenvalue.

   {\bf Proof of Theorem \ref{two-sided.theorem}.} 
   We only need to prove the upper bounds
    as  the lower bounds are from Theorem \ref{psi.lemma}. 
Let $s := |S|$ and let $u^*$ be defined by
 $$\phi_0^2 (L,S)  := s  \| \bigl < X_0 u^* \bigr > \|_2^2 . $$
 Then 
 $$ \hat \phi^2 (L,S) \le s  \| X u^* \|_{2,n}^2 =  \phi_0^2 (L,S) +s  \biggl ( \| X u^* \|_{2,n}^2
 - \| \bigl < X_0 u^* \bigr > \|_2^2 \biggr ) .$$
 But by Chebyshev's inequality, for all $t>0$
 $$ \PP \biggl  (  \|  X u^* \|_{2,n}^2
 - \| \bigl < X_0 u^* \bigr > \|_2^2 > \sqrt {t \over n}  \| \bigl < X_0 u^* \bigr > \|_4^2 \biggr ) \le 1/t . $$
 Insert the bound of Lemma \ref{fourthmoments.lemma} for $\| \bigl < X_0 u^* \bigr > \|_4/ \| \bigl < X_0 u^* \bigr > \|_2^4 $ or, in the case $m > 4$, 
 the bound
 $$\| \bigl < X_0 u^* \bigr > \|_4^4 \le C_m^{4} + \int_{C_m^{4}}^{\infty} \PP (
 | \bigl < X_0 u^* \bigr > | \ge t^{1/4} ) dt $$
 $$ \le C_m^{4} + C_m^m \int_{C_m^{4}}^{\infty}  t^{- m/4} dt 
  = C_m^{4} m/(m-4)  . $$
 This gives that with probability at least $1- 1/t$
 $$ s  \biggl ( \| X u^* \|_{2,n}^2
 - \| \bigl < X_0 u^* \bigr > \|_2^2 \biggr )  \le s \| \bigl < X_0 u \bigr >^* \|_2^2  \Delta_n^{\rm U} ( (L+1) \sqrt s, t)$$ $$=  \phi_0^2 (L,S)
 \Delta_n^{\rm U} ( (L+1) \sqrt s , t). $$
 The result for the restricted eigenvalue follows in the same way. 
 \hfill $\sqcup \mkern -12mu \sqcap$
 
 \subsection{Proofs for Section \ref{normalized.section}}
 
 We use the transfer principle to obtain lower bounds for
 sparse quadratic forms.
 
 {\bf Proof of Theorem \ref{transfer1.theorem}.} 
We apply result (\ref{first.equation}) of Theorem \ref{easy.theorem} to
$$ \inf \{ \| X u_J \|_{2,n}-1 : \ \| \bigl < X_0 u_J \bigr > \|_2 = 1 \} .$$
where 
$J$ is a fixed subset of $\{ 1 , \ldots , p \}$
with $|J| \le M^2$. There
are at most $p^{M^2}$ such subsets. Hence, by the union bound and replacing in 
the expression (\ref{first.equation}) of Theorem 
\ref{easy.theorem} the value $p$ by $M^2$ and $t$ by $t + M^2 \log p$
we have that with probability at least $1- \exp[-t]$
$$  \| X u_J \|_{2,n}^2 \ge \biggl (1- \bar \Delta_n^{\rm L} (M, t) \biggr ) \| \bigl < X_0 u_J \bigr >\|_2^2 \  \forall  \ u \in \R^p , \ \forall\ | J | \le 
M^2 . $$
The result follows now from Corollary \ref{transfer.corollary}.

\hfill $\sqcup \mkern -12mu \sqcap$

To handle the event ${\cal B}= \{ \max_j \hat \sigma_j^2 \le 1 + \epsilon \}$ we gave two lemmas.
Here are their proofs.

{\bf Proof of Lemma \ref{sigmaGauss.lemma}.}
Recall we assumed in the beginning of Section \ref{normalized.section} that $\| X_{0,j} \|_2 =1$ for all $j$.
 The assumption that the $X_{0,j}$ are sub-Gaussian implies
$$\| X_{0,j } \|_{\Psi_2} \le 2 C . $$ Hence
$$\EE | X_{0,j} |^{2k} \le k! (4C^2)^k  $$
and so 
$$ \EE | X_{0,j}^2 - \EE X_{0,j}^2 |^k \le 2^{k-1} k! (4C^2)^k= {k! \over 2} ( 8C^2)^k . $$
By Lemma 14.13 in \cite{BvdG2011} we find
$$ \PP \biggl ( \max_{1 \le j \le p} | \hat \sigma_j^2 - 1 |/ (8C^2)  \ge
 \sqrt {2 \log (2p) \over  n }  + \sqrt {2t \over n}+ { \log (2p) \over n }     + {t\over n } 
\biggr ) \le \exp [-t ] . $$
\hfill $\sqcup \mkern -12mu \sqcap$

{\bf Proof of Lemma \ref{sigmamoments.lemma}.} The moment conditions imply
 that
 $$\max_{1 \le j \le p } \| X_{0,j}^2  \|_k= \max_{1 \le j \le p } \| X_{0,j} \|_{2k}^2 \le
 \kappa_1^2 (2k)^{2 \alpha} = \kappa_1^2 2^{2 \alpha} k^{2 \alpha} , \ 1 \le k \le k_0/2  . $$
 But then
 $$ \max_{1 \le j \le p } \| X_{0,j}^2 -1 \|_k \le \kappa_1^2 2^{2 \alpha} k^{2 \alpha} +1 \le
 \kappa_1^2 2^{2 \alpha+1} k^{2 \alpha}, \ 1 \le k \le k_0/2 . $$
 We therefore have by Lemma \ref{Lecue.lemma}
 $$ \max_{1 \le j \le p } \|  \hat \sigma_j^2 - 1 \|_{k_0/2} \le c_0 \exp[ 4 \alpha -1] 
 \kappa_1^2 2^{2 \alpha+1} \sqrt {k_0 / (2n) } . $$
 It follows that
 $$ \biggl ( \EE \max_{1 \le j \le p } | \hat \sigma_j^2 -1 |^{k_0 /2} \biggr )^{2 \over  k_0} 
  \le p^{2 \over k_0 } c_0 \exp[ 4 \alpha -1] 
 \kappa_1^2 2^{2 \alpha+1} \sqrt {k_0 / (2n) }  $$ $$ = c_0 \exp[ 4 \alpha -1+ 2/\eta ] 
 \kappa_1^2 2^{2 \alpha+1} \sqrt {\eta \log p  / (2n) } := c_1 \sqrt {\log p / n }. $$
 But then by Chebyshev's inequality, for all $t >0$
 $$ \PP \biggl ( \max_{1 \le j \le p } \hat \sigma_j^2 - 1 \ge c_1 t^{2/k_0} \sqrt {\log p/n} \biggr  ) \le 
  1/t . $$
  
  \hfill $\sqcup \mkern -12mu \sqcap$
  
  Here is the proof for the case of normalized design.
  
  {\bf Proof of Theorem \ref{normalized.theorem}.} Consider the event
$${\cal A}:= \biggl \{  \| X u_J \|_{2,n}^2  \ge (1 - \Delta ) \| \bigl < X_0 u_J \bigr > \|_2^2  ,\
\forall \ u \in \R^p, \ J \subset \{ 1 , \ldots , p \} , \ | J | \le M^2 (\Delta)   \biggr \}  $$
Then on ${\cal A}$, by the transformation $u \mapsto \hat D^{-1/2} u $,
$$ u_J^T \hat R u_J  \ge (1 - \Delta ) u_J \hat D^{-1/2} \Sigma_0 \hat D^{-1/2} u_J  ,$$
for all $u \in \R^p, \ J \subset \{ 1 , \ldots , p \} , \ | J | \le M^2(\Delta) $. 
 The diagonal of 
 $$\hat D^{-1/2} ( \hat \Sigma - (1- \Delta) \Sigma_0) \hat D^{-1/2}$$
 is non-negative on ${\cal A}$ and less than or equal to $1 $.
 So on ${\cal A}$
by the transfer principle (Theorem \ref{transfer.theorem}) we know for all $u \in \R^p$
with $\| u \|_1^2 \le (L+1)^2  $
that
$$ u^T \hat R u \ge (1- \Delta) u^T \hat D^{-1/2} \Sigma_0 \hat D^{-1/2} u -  
{(L+1)^2  \over  M^2(\Delta) -1 }  . $$
We now note that
$$\inf_{\| u_S \|_1=1 , \ \| u_{-S } \|_1 \le L }
u^T \hat D^{-1/2} \Sigma_0 \hat D^{-1/2} u =
\inf_{\sum_{j \in S} \hat \sigma_j | u_j |= 1 , \ 
\sum_{j \notin S } \hat \sigma_j |u_j | \le L } u^T \Sigma_0 u . $$
But on ${\cal C}_S$
$$\sum_{j \in S} \hat \sigma_j |u_j | \le ( \sum_{j\in S}\hat \sigma_j^2 )^{1/2} \| u_S \|_2 \le
\sqrt {s (1+ \epsilon)} 
 \| u \|_2  . $$
Moreover on ${\cal A}$
$$\sum_{j \notin S } \hat \sigma_j |u_j | \ge (1- \sqrt \Delta) \| u_{-S} \|_1 . $$ 
Hence
$$\inf_{\| u_S \|_1=1 , \ \| u_{-S } \|_1 \le L }
u^T \hat D^{-1/2} \Sigma_0 \hat D^{-1/2} u \ge \inf_{ \| u_S \|_2 \ge 1/ \sqrt {s(1+ \epsilon)} , \
\| u_{-S} \|_1 \le L/ (1- \sqrt \Delta) } \| \bigl < X_0 u \bigr > \|_2^2 $$ 
$$ =  \inf_{ \| u_{S} \|_2 \ge 1 / \sqrt {s(1+ \epsilon)} , \
\| u_{-S} \|_1 \le L/ (1- \sqrt \Delta) }{  \| \bigl < X_0 u \bigr > \|_2^2 \over \| u_S \|_2^2 } \| u_S \|_2^2 $$
$$ \ge \inf_{\| u_S \|_2 =1 , \ \| u_{-S} \|_1 \le L\sqrt {s (1+ \epsilon)} / (1- \sqrt \Delta) } 
{\| \bigl < X_0 u \bigr > \|_2^2 \over s (1+ \epsilon) } =
{ \kappa_*^2 ( L \sqrt {1+ \epsilon} / (1- \sqrt \Delta) , S)  \over s (1+ \epsilon) } .  $$
The further bounds on the event ${\cal A} \cap {\cal B}_S$ follow in the same way.

\hfill $\sqcup \mkern -12mu \sqcap$

  \subsection{Proofs for Section \ref{almost-bounded.section}}
  
  We show that a vector $X_0$ which is $m$-th order strongly
  isotropic and has
 $p$ sub-Gaussian entries is  up to constants ``almost bounded" by $\sqrt {\log (2p)}$.
  
  {\bf Proof of Lemma \ref{truncation.lemma}.} The sub-Gaussianity implies that for all $j$ and all $s>0$
$$P ( |X_{0,j} | / C  \ge \sqrt {2s} ) \le 2 \exp[- s] . $$
We find 
$$P \biggl (  \max_{1 \le j \le p } | X_{0,j} | / C \ge \sqrt {2t+ 2 \log (2 p) +2m( \log n) / (m-2)  } \biggr )
$$ $$
\le 2p \exp [ -(t  + \log (2 p)  + 2m(\log n)  / (m-2)) ] = \exp[-t] n^{-{m \over m-2} }  . $$
The proof is finished by applying the inequality
$$ \| (\bigl < X_0 u \bigr > ){\rm l} \{ X_0 \notin A(t) \} \|_2^2 \le \tilde C_m^2 \biggl ( P ( X_0 \notin A(t)) \biggr )^{m-2 \over m}  . $$
\hfill $\sqcup \mkern -12mu \sqcap$

If we have $n$ independent copies of 
a vector $X_0$ which is $m$-th order strongly
  isotropic and has
 $p$ sub-Gaussian entries these $n \times p$ variables  are up to constants ``almost bounded" by $\sqrt {\log (np)}$.
 For such bounded random variables, we now prove to have uniform convergence of the
 empirical norm.

{\bf Proof of Theorem \ref{uniform.theorem}.}
Let $A := A(t)$ be defined as in Lemma \ref{truncation.lemma}.
Recall the notation
$$\| f(X) \|_{2,n}^2 := {1 \over n} \sum_{i=1}^n f^2(X_i) $$
so
$$ \| (Xu) {\rm l} \{ X \in A \} \|_{2,n}^2 := {1 \over n} \sum_{i=1}^n
(( X u)_i)^2 {\rm l} \{ X_i \in A\} . $$
Write
$$ \| X u \|_{2,n}^2 - \| \bigl < X_0 u \bigr > \|_2^2 =
\| (Xu) {\rm l} \{ X \in A \} \|_{2,n}^2 - \| (\bigl < X_0 u \bigr > ){\rm l} \{ X_0 \in A \} \|_2^2 
$$ $$ +
\| (Xu) {\rm l} \{ X \notin A \} \|_{2,n}^2 - \| (\bigl < X_0 u \bigr > ){\rm l} \{ X_0 \notin A \} \|_2^2
$$
We have
$$ \PP\biggl (  \| (Xu) {\rm l} \{ X \notin A \} \|_{2,n}^2 \not= 0 \biggr )  =
\PP ( \exists \ i: \ X_i \notin A ) $$
$$ = \PP ( \max_{1 \le i \le n } \max_{1 \le j \le p } | X_{i,j} | \ge
\sigma_X ( \sqrt {2t+ 2 \log (2p) + 2 m \log (n)/ (m-2)  } )$$ $$ \le n^{-{m \over m-2}}  \exp [-t] . $$
Moreover, by the method in 
\cite{vandeGeer2014},
for a universal
 constant $c_1$ with probability at least $1- \exp[-t] $
 $$ \sup_{\| \bigl < X_0 u \bigr > \|_2 \le 1 , \ \| u \|_1 \le M } 
 \biggl | \| (X u) {\rm l}   \{ X \in A \} \|_{2,n} - \| (\bigl < X_0 u \bigr >)  {\rm l} \{ X_0 \in A \} \|_2^2 \biggr | / c_1  $$ $$\le M \sigma_X 
  \sqrt { ( 2t +2 \log (2pn) ) (\log p \log^3 n + t)  \over n }
 $$ $$+ M^2 \sigma_X^2 { ( 2t +2 \log (2pn) )  ( \log p \log^3 n + t ) \over n } . $$
 \hfill $\sqcup \mkern -12mu \sqcap$
 
 \subsection{Proofs for Section \ref{martingale.section}}
 
 First comes the result for directed acyclic graphs.
 
 {\bf Proof of Lemma \ref{DAG.lemma}.}
We may write $X_0 = X_0 B + \epsilon_0  $ where $B= \{ \beta_{k,j} \}$ with
$\beta_{k,j} =0 $ for $k \ge j $.
Thus $X_0 (I-B)= \epsilon_0  $
so it suffices to show that $\epsilon_0 $ is sub-Gaussian with constant
$ C  $.
Note that for $k \not= j$, say $k<j$
$$\EE ( \epsilon_{0,j} \epsilon_{0,k} ) = \EE \biggl ( \epsilon_{0,k} \EE (\epsilon_{0,j} \vert {\cal F}_{j-1} ) \biggr ) = 0. $$
We let $\Omega^2 := {\rm diag} (\omega_1^2 , \ldots , \omega_p^2)$.
Suppose that $\|\Omega u \|_2= 1 $.
Then for all $\lambda \in \R$
$$ \EE \exp [ \lambda \epsilon_0 u  ] = \EE\exp [ \lambda \sum_{j=1}^p \epsilon_{0,j} u_j   ] \le
\exp[ \lambda^2 C^2 \sum_{j=1}^p u_j^2 \omega_j^2 / 2] $$ $$=
 \exp[ \lambda^2 \| \Omega u \|_2^2 C^2 / 2 ] 
= \exp[ \lambda^2 C^2/2 ] . $$
Hence, using ${\rm e}^{|x|} \le {\rm e}^x + {\rm e}^{-x}$,
for all $\lambda \ge 0$,
$$ \EE \exp [ \lambda| \epsilon_0 u |  ] \le 2 \exp[ \lambda^2 C^2/2 ] . $$

\hfill $\sqcup \mkern -12mu \sqcap$

We now prove isotropy under conditional sub-Gaussian assumptions.

{\bf Proof of Lemma \ref{sub-Gauss.lemma}.}
We clearly have for all $\lambda \in \R $ and all $u \in \R^p$
$$\EE\biggl  ( \exp\biggl [ \lambda \bigl < X_0 u \bigr >  - \lambda^2 \sum_{j=1}^p { u_j^2 V_j^2 } / 2 
\biggr ]\biggr  \vert {\cal F}_0 \biggr  ) \le 1 . $$
If the $\{ V_j \}_{j=1}^p $ are ${\cal F}_0$-measurable this gives
$$\EE ( \exp [ \lambda \bigl < X_0 u \bigr >  ] \vert {\cal F}_0 ) \le \exp \biggl [ \lambda^2 \sum_{j=1}^p  u_j^2 V_j^2  / 2  \biggr ] . $$
We now use that (see e.g. \cite{BvdG2011}, Lemma 14.7)
$$\EE(  | \bigl < X_0 u \bigr > |^m \vert {\cal F}_0 )  \le \biggl ( { m \over \lambda}  +
{\lambda \sum_{j=1}^p u_j^2 V_j^2 \over 2 } \biggr )^m , $$
and we choose $\lambda = \sqrt {2m} /( \sum_{j=1}^p u_j V_j )^{1/2}$.
This gives
$$\EE ( | \bigl < X_0 u \bigr > |^m \vert {\cal F}_0 ) \le ( 2m)^{m/2} ( \sum_{j=1}^p u_j V_j )^{m/2} . $$
But then
$$\EE | \bigl < X_0 u \bigr > |^m \le ( 2m)^{m/2} \| u \|_2^m \mu_m^m . $$

For the case where $\{ V_j \}_{j=1}^p $ is predictable, we use that 
$$\EE\biggl  ( \exp\biggl [ \lambda \bigl < X_0 u \bigr >  - \lambda^2 \sum_{j=1}^p { u_j^2 V_j^2 } / 2 
\biggr ] \biggr  ) \le 1 . $$
and hence by standard arguments, for any positive
$a$ and $b$
$$ \PP \biggl ( { |\bigl < X_0 u \bigr >| \over \| u \|_2 \mu_m } \ge a , \  { \sum_{j=1}^p u_j^2 V_j^2 \over
\| u \|_2^2 \mu_m^2 } \le b^2 \biggr ) \le 2 \exp [ - a^2 / (2b^2) ] . $$
Choosing $a= b \sqrt {2m \log b}$ and $b= {\rm e}^{s/m} $ gives
$$ \PP \biggl ( { |\bigl < X_0 u \bigr >| \over \| u \|_2 \mu_m } \ge {\rm e}^{s/m} \sqrt {2s}  , \  { \sum_{j=1}^p u_j^2 V_j^2 \over
\| u \|_2^2 \mu_m^2 } \le {\rm e}^{2s/m}  \biggr ) \le 2 {\rm e}^{-s}  . $$
We thus find
$$ \PP \biggl ( { |\bigl < X_0 u \bigr >| \over \| u \|_2 \mu_m } \ge {\rm e}^{s/m} \sqrt {2s} \biggr ) \le 3 {\rm e}^{-s} . $$
We have
$$ { \EE | \bigl < X_0 u \bigr > |^{m_0}  \over \| u \|_2^{m_0}  \mu_m^{m_0} } =
\int_0^{\infty}  \PP \biggl ( { |\bigl < X_0 u \bigr >| \over \| u \|_2 \mu_m } > t^{1/ m_0} \biggr ) dt$$
$$ = \int_0^{\infty}  \PP \biggl ( { |\bigl < X_0 u \bigr >| \over \| u \|_2 \mu_m } > {\rm e}^{s/m} \sqrt {2s} \biggr ) 
d ({\rm e}^{s/m} \sqrt {2s} )^{m_0}$$
$$ \le  3  \int_0^{\infty} {\rm e}^{-s} 
d ({\rm e}^{s/m} \sqrt {2s} )^{m_0}$$
$$ = 3 \int_0^{\infty} ({\rm e}^{s/m} \sqrt {2s} )^{m_0} {\rm e}^{-s} ds $$ $$ =
3 (2^{m_0/2})  \biggl ( { m \over m- m_0 }\biggr ) ^{m_0/2 +1}   \Gamma (m_0 / 2 +1) . $$
\hfill $ \sqcup \mkern -12mu \sqcap$

The final proof concerns isotropy under conditional sub-exponential assumptions.

{\bf Proof of Lemma \ref{Bernstein.lemma}.}
We invoke the inequality
$$\log \EE \exp [Z] \le \EE {\rm e}^{|Z|} -1- \EE |Z| $$
(see e.g. Lemma 14.1 in \cite{BvdG2011}) which holds
for a random variable $Z \in \R$ with mean zero. 
Moreover
$$\EE {\rm e}^{|Z|} -1- \EE |Z| = \sum_{k=2}^{\infty} \EE |Z|^k / k! . $$
By the Bernstein condition one readily sees that
for all $\lambda\in \R $ and $u \in \R^p$ with $ |\lambda |  K \| u \|_2  < 1$
$$\EE \biggl ( \exp \biggl [ \lambda  \bigl < X_0 u \bigr >  - \lambda^2 \sum_{j=1}^p { u_j^2 V_j^2 }  / (2(1- |\lambda| K \| u \|_2 )  
\biggr ] \biggr \vert {\cal F}_0 \biggr ) \le 1 . $$
If the $\{ V_j \}_{j=1}^p $ are ${\cal F}_0$-measurable this gives
$$\EE ( \exp[ \lambda | \bigl < X_0 u \bigr > | ] \vert {\cal F}_0 ) \le \exp \biggl [ \lambda^2 \sum_{j=1}^p { u_j^2 V_j^2 }  / (2(1- |\lambda| K \| u \|_2 )  \biggr ] . $$
So then (see e.g. \cite{BvdG2011}, Lemma 14.7) for $0 < \lambda < K \| u \|_2 $
$$ \EE ( | \bigl < X_0 u \bigr > |^m \vert {\cal F}_0 )  \le \biggl ({ m \over \lambda }+ { \lambda \sum_{j=1}^p u_j^2 V_j^2 \over
2 ( 1- \lambda K \| u \|_2 ) } \biggr )^m . $$
Now choose 
$$ {1 \over \lambda } = K \| u \|_2 + \biggl ( {\sum_{j=1}^p u_j^2 V_j^2 \over 2 m } \biggr )^{1/2} . $$
Then we get
$$ \EE ( | \bigl < X_0 u \bigr > |^m \vert {\cal F}_0 )  \le 
\biggl ( \sqrt {2m} (\sum_{j=1}^p u_j^2 V_j^2 )^{1/2}  + m \| u \|_2 K  \biggr )^m .$$
This implies the result for non-random  $\{ V_j \}_{j=1}^p$.
If they are ${\cal F}_0$-measurable we find 
$$\EE  | \bigl < X_0 u \bigr > |^m \le 2^{m-1} \biggl ( \sqrt {2m} \| u \|_2 \mu_m \biggr )^m  +2^{m-1}  \biggl ( m \| u \|_2 K \biggr  )^m . $$

If the $\{ V_j \}_{j=1}^p$ are only predictable, we use that for all positive $a$ and $b$
and for $\tilde K := K/ \mu_m$
$$\PP \biggl ({ | \bigl < X_0 u \bigr > | \over \| u \|_2 \mu_m } \ge b \sqrt {2a} + \tilde K a  , \ 
{\sum_{j=1}^p u_j^2 V_j^2 \over \| u \|_2^2 \mu_m^2 } \le b^2  \biggr ) \le 2 \exp [-a] . $$
Write $a=s $ and $b = {\rm e}^{s/m} $ to find that
$$\PP \biggl ({ | \bigl < X_0 u \bigr > | \over \| u \|_2 \mu_m } \ge {\rm e}^{s/m}  \sqrt {2s} + \tilde K s  , \ 
{\sum_{j=1}^p u_j^2 V_j^2 \over \| u \|_2^2 \mu_m^2 } \le {\rm e}^{2s/m}  \biggr ) \le 2 \exp[-s]  $$
and so
$$\PP \biggl ({ | \bigl < X_0 u \bigr > | \over \| u \|_2 \mu_m } \ge {\rm e}^{s/m}  \sqrt {2s} + \tilde K s   \biggr ) \le 3 \exp[-s] . $$

It follows that
$$ \EE \biggl ({  | \bigl < X_0 u \bigr > | \over \| u \|_2 \mu_m } \biggr )^{m_0} =
 \int_0^{\infty} \PP \biggl ( { | \bigl < X_0 u \bigr > | \over \| u \|_2 \mu_m } \ge t^{1 / m_0} \biggr ) dt $$
 $$ \le 3 \int_0^{\infty} \biggr ( {\rm e}^{s/m}  \sqrt {2s} + \tilde K s   \biggr )^{m_0} {\rm e}^{-s} ds .$$
 But
$$  \int_0^{\infty} \biggr ( {\rm e}^{s/m}  \sqrt {2s}   \biggr )^{m_0} {\rm e}^{-s} ds=
  (2^{m_0/2} ) \biggl ( { m \over m- m_0 } \biggr )^{m_0 /2 +1} \Gamma
 (m_0/2 +1 ) $$
 and
 $$\int_0^{\infty} \biggr (   \tilde K s   \biggr )^{m_0} {\rm e}^{-s} ds= \tilde K^{m_0} 
 \Gamma (m_0 + 1) . $$
 Hence by the triangle inequality
 $$ { \|  \bigl < X_0 u \bigr >  \|_{m_0} \over \| u \|_2 \mu_m } \le 
  \sqrt { 2m \over m- m_0 } \biggl ( { 3 m \Gamma (m_0 /2 + 1)  \over m-m_0}  \biggr )^{m_0 /2 +1}+
 \tilde K 
\biggr (3  \Gamma (m_0 + 1) \biggr )^{1/ m_0} . $$

\hfill $\sqcup \mkern -12mu \sqcap$

\bibliographystyle{plainnat}
\bibliography{reference}

\end{document}